\documentclass[reqno, letter, 10pt]{amsart} 

\usepackage[english]{babel}
\usepackage[utf8x]{inputenc}
\usepackage{mathtools}
\mathtoolsset{showonlyrefs,showmanualtags}
\usepackage{empheq}
\usepackage{amsmath, amssymb}
\usepackage{stmaryrd}
\usepackage{enumerate,enumitem}
\usepackage[perpage,multiple,bottom]{footmisc}
\usepackage[usenames,dvipsnames]{color}
\usepackage{hyperref}
\hypersetup{
	colorlinks=true,
	citecolor=blue,
	linkcolor=blue,
	filecolor=blue, 
	urlcolor=blue,
}
\usepackage[toc,page]{appendix}
\usepackage{xcolor}
\usepackage[breakable]{tcolorbox}
\usepackage{tikz}
\usetikzlibrary{shapes, backgrounds} 
\usepackage{tikz-cd}
\definecolor{orange}{rgb}{1,0.5,0}
\usepackage[all,cmtip,knot]{xy}

\newcommand{\nc}{\newcommand}
\newcommand{\rnc}{\renewcommand}

\nc{\andreacomment}[1]{\footnote{\color{red}{Andrea: #1}}}
\nc{\valeriocomment}[1]{\footnote{\color{blue}{Valerio: #1}}}
\nc{\omitvaleriocomment}[1]{}
\nc{\Omit}[1]{}
\nc{\tocheck}[1]{{\color{magenta} #1}}

\nc{\triend}{\parbox{2mm}{\hfill} \hfill\text{\hspace{0.2mm}}\hfill$\triangle$}
\nc{\ocend}{\parbox{2mm}{\hfill} \hfill\text{\hspace{0.2mm}}\hfill$\oslash$}
\nc{\black}{\parbox{2mm}{\hfill}\text{\hspace{0.2mm}}\hfill$\blacksquare$}

\numberwithin{equation}{section}

\newtheorem*{theorem}{Theorem}
\newtheorem*{proposition}{Proposition}
\newtheorem*{lemma}{Lemma}
\newtheorem*{corollary}{Corollary}

\newenvironment{definition}{\medskip\paragraph{{\bf Definition.}}}{\par\medskip} 
\newenvironment{remark}{\medskip\paragraph{{\bf Remark.}}}{\par\medskip}
\newenvironment{remarks}{\medskip\paragraph{{\bf Remarks.} }}{\par\medskip}
\nc{\cA}{{\mathcal A}}
\nc{\cB}{{\mathcal B}}
\nc{\cC}{{\mathcal C}}
\nc{\cD}{{\mathcal D}}
\nc{\cE}{{\mathcal E}}
\nc{\cF}{{\mathcal F}}
\nc{\cG}{{\mathcal G}}
\nc{\cH}{{\mathcal H}}
\nc{\cI}{{\mathcal I}}
\nc{\cJ}{{\mathcal J}}
\nc{\cK}{{\mathcal K}}
\nc{\cL}{{\mathcal L}}
\nc{\cM}{\mathcal{M}}
\nc{\cN}{{\mathcal N}}
\nc{\cO}{{\mathcal O}}
\nc{\cP}{{\mathcal P}}
\nc{\cQ}{{\mathcal Q}}
\nc{\cR}{{\mathcal R}}
\nc{\cS}{\mathcal{S}}
\nc{\cT}{{\mathcal T}}
\nc{\cU}{\mathcal{U}}
\nc{\cV}{{\mathcal V}}
\nc{\cX}{{\mathcal X}}
\nc{\cY}{\mathcal{Y}}
\nc{\cW}{\mathcal{W}}
\nc{\cZ}{{\mathcal Z}}

\nc{\bbA}{{\mathbb{A}}}
\nc{\bbB}{{\mathbb{B}}}
\nc{\bbC}{{\mathbb{C}}}
\nc{\bbD}{{\mathbb{D}}}
\nc{\bbE}{{\mathbb{E}}}
\nc{\bbF}{{\mathbb{F}}}
\nc{\bbG}{{\mathbb{G}}}
\nc{\bbH}{{\mathbb{H}}}
\nc{\bbI}{{\mathbb{I}}}
\nc{\bbJ}{{\mathbb{J}}}
\nc{\bbK}{{\mathbb{K}}}
\nc{\bbL}{{\mathbb{L}}}
\nc{\bbM}{{\mathbb{M}}}
\nc{\bbN}{{\mathbb{N}}}
\nc{\bbO}{{\mathbb{O}}}
\nc{\bbP}{{\mathbb{P}}}
\nc{\bbQ}{{\mathbb{Q}}}
\nc{\bbR}{{\mathbb{R}}}
\nc{\bbS}{{\mathbb{S}}}
\nc{\bbT}{{\mathbb{T}}}
\nc{\bbU}{{\mathbb{U}}}
\nc{\bbV}{{\mathbb{V}}}
\nc{\bbX}{{\mathbb{X}}}
\nc{\bbY}{{\mathbb{Y}}}
\nc{\bbW}{{\mathbb{W}}}
\nc{\bbZ}{{\mathbb{Z}}}

\usepackage{mathrsfs}

\nc{\scrA}{{\mathscr A}}
\nc{\scrB}{{\mathscr B}}
\nc{\scrC}{{\mathscr C}}
\nc{\scrD}{{\mathscr D}}
\nc{\scrE}{{\mathscr E}}
\nc{\scrF}{{\mathscr F}}
\nc{\scrG}{{\mathscr G}}
\nc{\scrH}{{\mathscr H}}
\nc{\scrI}{{\mathscr I}}
\nc{\scrJ}{{\mathscr J}}
\nc{\scrK}{{\mathscr K}}
\nc{\scrL}{{\mathscr L}}
\nc{\scrM}{{\mathscr M}}
\nc{\scrN}{{\mathscr N}}
\nc{\scrO}{{\mathscr O}}
\nc{\scrP}{{\mathscr P}}
\nc{\scrQ}{{\mathscr Q}}
\nc{\scrR}{{\mathscr R}}
\nc{\scrS}{{\mathscr S}}
\nc{\scrT}{{\mathscr T}}
\nc{\scrU}{{\mathscr U}}
\nc{\scrV}{{\mathscr V}}
\nc{\scrX}{{\mathscr X}}
\nc{\scrY}{{\mathscr Y}}
\nc{\scrW}{{\mathscr W}}
\nc{\scrZ}{{\mathscr Z}}

\nc{\sfA}{{\mathsf A}}
\nc{\sfB}{{\mathsf B}}
\nc{\sfC}{{\mathsf C}}
\nc{\sfD}{{\mathsf D}}
\nc{\sfE}{{\mathsf E}}
\nc{\sfF}{{\mathsf F}}
\nc{\sfG}{{\mathsf G}}
\nc{\sfH}{{\mathsf H}}
\nc{\sfI}{{\mathsf I}}
\nc{\sfJ}{{\mathsf J}}
\nc{\sfK}{{\mathsf K}}
\nc{\sfL}{{\mathsf L}}
\nc{\sfM}{{\mathsf M}}
\nc{\sfN}{{\mathsf N}}
\nc{\sfO}{{\mathsf O}}
\nc{\sfP}{{\mathsf P}}
\nc{\sfQ}{{\mathsf Q}}
\nc{\sfR}{{\mathsf R}}
\nc{\sfS}{{\mathsf S}}
\nc{\sfT}{{\mathsf T}}
\nc{\sfU}{{\mathsf U}}
\nc{\sfV}{{\mathsf V}}
\nc{\sfX}{{\mathsf X}}
\nc{\sfY}{{\mathsf Y}}
\nc{\sfW}{{\mathsf W}}
\nc{\sfZ}{{\mathsf Z}}

\nc{\sfa}{{\mathsf a}}
\nc{\sfb}{{\mathsf b}}
\nc{\sfc}{{\mathsf c}}
\nc{\sfd}{{\mathsf d}}
\nc{\sfe}{{\mathsf e}}
\nc{\sff}{{\mathsf f}}
\nc{\sfg}{{\mathsf g}}
\nc{\sfh}{{\mathsf h}}
\nc{\sfi}{{\mathsf i}}
\nc{\sfj}{{\mathsf j}}
\nc{\sfk}{{\mathsf k}}
\nc{\sfl}{{\mathsf l}}
\nc{\sfm}{{\mathsf m}}
\nc{\sfn}{{\mathsf n}}
\nc{\sfo}{{\mathsf o}}
\nc{\sfp}{{\mathsf p}}
\nc{\sfq}{{\mathsf q}}
\nc{\sfr}{{\mathsf r}}
\nc{\sfs}{{\mathsf s}}
\nc{\sft}{{\mathsf t}}
\nc{\sfu}{{\mathsf u}}
\nc{\sfv}{{\mathsf v}}
\nc{\sfx}{{\mathsf x}}
\nc{\sfy}{{\mathsf y}}
\nc{\sfw}{{\mathsf w}}
\nc{\sfz}{{\mathsf z}}

\nc {\bfA}{{\mathbf A}}
\nc {\bfB}{{\mathbf B}}
\nc {\bfC}{{\mathbf C}}
\nc {\bfD}{{\mathbf D}}
\nc {\bfE}{{\mathbf E}}
\nc {\bfF}{{\mathbf F}}
\nc {\bfG}{{\mathbf G}}
\nc {\bfH}{{\mathbf H}}
\nc {\bfI}{{\mathbf I}}
\nc {\bfJ}{{\mathbf J}}
\nc {\bfK}{{\mathbf K}}
\nc {\bfL}{{\mathbf L}}
\nc {\bfM}{{\mathbf M}}
\nc {\bfN}{{\mathbf N}}
\nc{\bfO}{{\mathbf O}}
\nc {\bfP}{{\mathbf P}}
\nc {\bfQ}{{\mathbf Q}}
\nc {\bfR}{{\mathbf R}}
\nc {\bfS}{{\mathbf S}}
\nc {\bfT}{{\mathbf T}}
\nc {\bfU}{{\mathbf U}}
\nc {\bfV}{{\mathbf V}}
\nc {\bfX}{{\mathbf X}}
\nc {\bfY}{{\mathbf Y}}
\nc {\bfW}{{\mathbf W}}
\nc {\bfZ}{{\mathbf Z}}

\nc {\fka}{{\mathfrak a}}
\nc {\fkb}{{\mathfrak b}}
\nc {\fkc}{{\mathfrak c}}
\nc {\fkd}{{\mathfrak d}}
\nc {\fke}{{\mathfrak e}}
\nc {\fkf}{{\mathfrak f}}
\nc {\fkg}{{\mathfrak g}}
\nc {\fkh}{{\mathfrak h}}
\nc {\fki}{{\mathfrak i}}
\nc {\fkj}{{\mathfrak j}}
\nc {\fkk}{{\mathfrak k}}
\nc {\fkl}{{\mathfrak l}}
\nc {\fkm}{{\mathfrak m}}
\nc {\fkn}{{\mathfrak n}}
\nc {\fko}{{\mathfrak o}}
\nc {\fkp}{{\mathfrak p}}
\nc {\fkq}{{\mathfrak q}}
\nc {\fkr}{{\mathfrak r}}
\nc {\fks}{{\mathfrak s}}
\nc {\fkt}{{\mathfrak t}}
\nc {\fku}{{\mathfrak u}}
\nc {\fkv}{{\mathfrak v}}
\nc {\fkx}{{\mathfrak x}}
\nc {\fky}{{\mathfrak y}}
\nc {\fkw}{{\mathfrak w}}
\nc {\fkz}{{\mathfrak z}}

\nc {\fkA}{{\mathfrak A}}
\nc {\fkB}{{\mathfrak B}}
\nc {\fkC}{{\mathfrak C}}
\nc {\fkD}{{\mathfrak D}}
\nc {\fkE}{{\mathfrak E}}
\nc {\fkF}{{\mathfrak F}}
\nc {\fkG}{{\mathfrak G}}
\nc {\fkH}{{\mathfrak H}}
\nc {\fkI}{{\mathfrak I}}
\nc {\fkJ}{{\mathfrak J}}
\nc {\fkK}{{\mathfrak K}}
\nc {\fkL}{{\mathfrak L}}
\nc {\fkM}{{\mathfrak M}}
\nc {\fkN}{{\mathfrak N}}
\nc {\fkO}{{\mathfrak O}}
\nc {\fkP}{{\mathfrak P}}
\nc {\fkQ}{{\mathfrak Q}}
\nc {\fkR}{{\mathfrak R}}
\nc {\fkS}{{\mathfrak S}}
\nc {\fkT}{{\mathfrak T}}
\nc {\fkU}{{\mathfrak U}}
\nc {\fkV}{{\mathfrak V}}
\nc {\fkX}{{\mathfrak X}}
\nc {\fkY}{{\mathfrak Y}}
\nc {\fkW}{{\mathfrak W}}
\nc {\fkZ}{{\mathfrak Z}}

\nc{\Hom}{\operatorname{Hom}}
\nc{\Ext}{\operatorname{Ext}}
\nc{\Aut}{\operatorname{Aut}}
\nc{\End}{\operatorname{End}}
\nc{\Tor}{\operatorname{Tor}}
\nc{\Pic}{\operatorname{Pic}}
\nc{\Ind}{\operatorname{Ind}}
\nc{\Fun}{\operatorname{Fun}}
\nc{\ev}{\operatorname{ev}}
\nc{\Res}{\operatorname{Res}}
\nc{\Obj}{\operatorname{Obj}}
\nc{\Rep}{\operatorname{Rep}}
\nc{\Ad}{\operatorname{Ad}}
\nc{\Ker}{\operatorname{Ker}}

\nc{\Repss}{\operatorname{Rep}^{\scsop{ss}}}

\nc{\id}{{\operatorname{id}}}
\nc{\Spec}{\operatorname{Spec}}
\nc{\fib}{\operatorname{fib}}
\nc{\colim}{\operatorname{colim}}
\nc{\rank}{\operatorname{rk}}
\nc{\opspan}{\operatorname{span}}

\nc {\ul}{\underline}
\nc {\ol}{\overline}
\nc {\wt}{\widetilde}
\nc {\wh}{\widehat}
\nc{\scs}{\scriptscriptstyle}
\nc{\scsop}{\scriptscriptstyle\operatorname} 
\nc{\op}{{\opertorname{op}}}

\nc {\ie}{{\emph{i.e.}}, }
\nc {\eg}{{\emph{e.g.}}, }
\nc {\aand}{\qquad\mbox{and}\qquad}

\nc{\drc}[1]{\delta_{#1}}
\nc{\der}{\partial}
\nc{\sfad}{\mathsf{ad}}
\nc{\ten}{\otimes}

\rnc{\a}{\fka}
\rnc{\b}{\fkb}
\rnc{\c}{\fkc}
\rnc{\d}{\fkd}
\nc{\e}{\fke}
\nc{\f}{\fkf}
\nc{\g}{\fkg}
\nc{\h}{\fkh}
\rnc{\i}{\fki}
\rnc{\j}{\fkj}
\rnc{\k}{\fkk}
\rnc{\l}{\fkl}
\nc{\m}{\fkm}
\nc{\n}{\fkn}
\rnc{\o}{\fko}
\nc{\p}{\fkp}
\nc{\q}{\fkq}
\rnc{\r}{\fkr}
\nc{\s}{\fks}
\rnc{\t}{\fkt}
\rnc{\u}{\fku}
\rnc{\v}{\fkv}
\nc{\x}{\fkx}
\nc{\y}{\fky}
\nc{\w}{\fkw}
\nc{\z}{\fkz}

\nc{\A}{\cA}
\nc{\B}{\cB}
\nc{\C}{\cC}
\nc{\D}{\cD}
\nc{\E}{\cE}
\nc{\F}{\cF}
\nc{\G}{\cG}
\rnc{\H}{\cH}
\nc{\I}{\cI}
\nc{\J}{\cJ}
\nc{\K}{\cK}
\rnc{\L}{\cL}
\nc{\M}{\cM}
\nc{\N}{\cN}
\rnc{\O}{\cO}
\rnc{\P}{\cP}
\nc{\Q}{\cQ}
\nc{\R}{\cR}
\rnc{\S}{\cS}
\nc{\T}{\cT}
\nc{\U}{\cU}
\nc{\V}{\cV}
\nc{\X}{\cX}
\nc{\Y}{\cY}
\nc{\W}{\cW}
\nc{\Z}{\cZ}

\nc{\IC}{\bbC}
\nc{\IR}{\bbR}
\nc{\IQ}{\bbQ}
\nc{\IZ}{\bbZ}
\nc{\IV}{\bbV}

\nc{\nablah}{\sfh}

\nc{\hext}[1]{{#1}[\negthinspace[\hbar]\negthinspace]}

\nc{\dgr}{\bbD} 
\nc{\vect}{\operatorname{Vect}}

\nc{\GCM}[1]{{\mathsf #1}} 
\nc{\Diag}{\operatorname{diag}}
\nc{\symd}[1]{d_{#1}}
\nc{\symdi}[1]{d^{-1}_{#1}}
\nc{\ess}{^{\operatorname{e}}}

\nc{\Par}{\operatorname{Par}} 
\nc{\Ql}{\mathsf{Q}}
\nc{\Qp}{\mathsf{Q}_{+}}
\nc{\Qv}{\mathsf{Q}^\vee}

\nc{\Pl}{\mathsf{P}}
\nc{\Pp}{\mathsf{P}_+}

\nc{\hgt}{\operatorname{ht}} 
\nc{\rsm}[1]{\operatorname{m}_{#1}}

\nc{\Rs}[1]{\Delta_{#1}} 

\nc{\iip}[2]{\langle#1,#2\rangle}
\nc{\mult}{\operatorname{mult}}

\nc{\rt}[1]{\alpha_{#1}} 
\nc{\rtv}[2]{e_{#1}^{#2}} 
\nc{\crt}[1]{h_{#1}} 
\nc{\tinv}[1]{t_{#1}} 

\nc{\rootsys}{\mathsf{\Delta}}
\nc{\Rp}{\Delta_+}

\nc{\bpm}[1]{\b^{\pm}_{#1}}
\nc{\bmp}[1]{\b^{\mp}_{#1}}
\nc{\bp}[1]{\b^{+}_{#1}}
\nc{\bm}[1]{\b^{-}_{#1}}
\nc{\npm}[1]{\n^{\pm}_{#1}}
\nc{\nmp}[1]{\n^{\mp}_{#1}}
\nc{\np}[1]{\n^{+}_{#1}}
\nc{\nm}[1]{\n^{-}_{#1}}

\nc{\Uh}{U_{\hbar}}
\nc{\Uhg}{\Uh\g}
\nc{\Uqg}{U_q\g}
\nc{\Uhbpm}[1]{\Uh\bpm{#1}}
\nc{\Uhbmp}[1]{\Uh\bmp{#1}}
\nc{\Uhbp}[1]{\Uh\bp{#1}}
\nc{\Uhbm}[1]{\Uh\bm{#1}}
\nc{\Uqbm}[1]{U_q\bm{#1}}
\nc{\Uhnpm}[1]{\Uh\npm{#1}}
\nc{\Uhnmp}[1]{\Uh\nmp{#1}}
\nc{\Uhnp}[1]{\Uh\np{#1}}
\nc{\Uqnp}[1]{U_q\np{#1}}
\nc{\Uhnm}[1]{\Uh\nm{#1}}

\nc{\qKi}[2]{q_{#1}^{{#2}\crt{#1}}} 

\nc{\Rmx}{\mathbf{R}} 
\nc{\qWSk}[1]{{\mathbf{S}}_{#1}} 
\nc{\qWS}[1]{\mathbf{S}_{\C,#1}} 
\nc{\qWT}[1]{\mathbf{T}_{#1}} 

\nc{\texp}[1]{\wt{s}_{#1}}
\nc{\texpb}[1]{\tau(#1)}
\nc{\topS}[1]{S_{#1}} 

\nc{\Br}[1]{\mathcal{B}_{#1}} 
\nc{\PBr}[1]{\mathcal{P}_{#1}} 
\nc{\ulm}{\ul{m}}
\nc{\BDm}{\Br{\dgr}^{\ulm}}

\nc {\KM}{Kac--Moody }
\nc {\KMA}{Kac--Moody algebra }
\nc {\Bw}{\Br{W}}
\nc {\Pw}{\PBr{W}}
\nc {\PwJ}{\PBr{W_\bfJ}}
\renewcommand {\int}{^{\scriptscriptstyle{\operatorname{int}}}}
\nc {\sfointh}{\mathsf{f}\int_\hbar}
\nc {\sffOh}{\mathsf{f}_\hbar}
\nc{\Oh}{\O_\hbar}
\nc {\veps}{\varepsilon}
\nc {\sfPhi}{\mathsf{\Phi}}
\nc {\vech}{\mathsf{Vec}_\hbar}
\nc {\re}{^{\scriptscriptstyle{\operatorname{re}}}}
\nc {\sint}{\scsop{int}}

\rnc{\sl}[1]{\mathfrak{sl}_{#1}}

\nc{\Cu}[1]{\mathcal{C}_{#1}}
\nc{\ku}[1]{\mathcal{K}_{#1}}
\nc{\Ku}[2]{\mathcal{K}_{#1}^{#2}}

\nc{\pp}{}

\nc{\cor}[1]{\crt{#1}}
\nc{\Oinf}{\O_{\infty}}
\nc{\Oinfp}{\O_{\infty}^{\pp}}
\nc{\Oint}{\O\int}
\nc{\Ointinf}{\O_{\infty}^{\sint}}
\nc{\Oinfint}{\O_{\infty}^{\scsop{int}}}
\nc{\hOinf}{\O_{\infty}}
\nc{\hOinfint}{\O_{\infty}^{\scsop{int}}}

\nc{\Og}{\O_{\g}}
\nc{\Oinfg}{\O_{\infty,\g}}
\nc{\Ointg}{\O\int_{\g}}
\nc{\Oinfintg}{\O_{\infty,\g}^{\sint}}
\nc{\OUhg}{\O_{\Uhg}}
\nc{\OinfUhg}{\O_{\infty,\Uhg}}
\nc{\OinfUhgp}{\O^{\pp}_{\infty,\Uhg}}
\nc{\OinfUhgD}[1]{\O_{\infty,\Uhg_{#1}}}
\nc{\OinfintUhg}{\O_{\infty,\Uhg}^{\scsop{int}}}
\nc{\OinfintUhgD}[1]{\O_{\infty,\Uhg_{#1}}^{\sint}}

\nc{\Ohg}{\O^{\hbar}_{\g}}
\nc{\Ohintg}{\O^{\hbar,\sint}_{\g}}
\nc{\Ohinfintg}{\O_{\infty,\g}^{\hbar,\sint}}
\nc{\Ohinfg}{\O_{\infty,\g}^{\hbar}}
\nc{\Ohinfgp}{\O_{\infty,\g}^{\hbar,\pp}}

\nc{\OhinfintgD}[1]{\O_{\infty,\g_{#1}}^{\hbar,\sint}}
\nc{\OhinfgD}[1]{\O_{\infty,\g_{#1}}^{\hbar}}

\nc{\OinfintUhsli}{\O_{\infty,U_\hbar\mathfrak{sl}_2^{\alpha_i}}^{\scsop{int}}}

\nc{\CQUOinf}{\End(\FF{\hbar})}
\nc{\CQUOinfx}{\Aut(\FF{\hbar})}
\nc{\CQUOinfint}{\End(\FF{\hbar}^{\sint})}
\nc{\CQUOinfintx}{\Aut(\FF{\hbar}^{\sint})}

\nc{\CQUOinfintD}[1]{\End(\FF{\hbar, #1}^{\sint})}
\nc{\CQUOinfintxD}[1]{\Aut(\FF{\hbar, #1}^{\sint})}

\nc{\CUOhinf}{\End(\FF{})}
\nc{\CUOhinfx}{\Aut(\FF{})}
\nc{\CUOhinfint}{\End(\FF{}^{\sint})}
\nc{\CUOhinfintx}{\Aut(\FF{}^{\sint})}
\nc{\CUOhinfD}[1]{\End(\FF{#1})}
\nc{\CUOhinfxD}[1]{\Aut(\FF{#1})}
\nc{\CUOhinfintD}[1]{\End(\FF{#1}^{\sint})}
\nc{\CUOhinfintxD}[1]{\Aut(\FF{#1})^{\sint}}

\nc{\CUOhinfintBD}[2]{\End(\FF{#1#2}^{\sint})}

\nc{\CUOhss}{\End(\FF{\h})}

\nc {\Dg}{D}
\nc {\DrA}{\D}
\nc {\DrAh}{{\mathcal D}_\hbar}
\nc {\DrAq}{{\mathcal D}_q}
\nc {\DrQh}{{\mathcal Q}_\hbar}
\nc {\DrAhP}{{\mathcal D}_\hbar^{\Pl}}
\nc {\Vect}{\operatorname{Vect}}
\nc {\Vecth}{\Vect_{\hbar}}

\nc{\Dr}{\D}

\nc {\kalpha}{\kappa_{\alpha}}
\nc {\kbeta}{\kappa_{\beta}}
\nc {\calpha}{C_{\alpha}}
\nc {\calkalpha}{{\mathcal K}_\alpha}
\nc {\nablac}{\nabla_{C}}
\nc {\nablak}{\nabla_{\K}}
\nc {\wtnablak}{\wt{\nabla}_{\K}}
\nc {\nablawk}{\nabla_{:\K:}}

\nc{\VA}{V_\A}

\nc {\PT}{\mathscr{P}}
\nc {\bfa}{\mathbf{a}}
\nc {\Poid}{\mathbf{\Pi}_1}
\nc {\WPoid}[1]{W\ltimes\Poid{#1}}

\nc{\aw}{\mathscr{A }}
\nc{\bw}{\mathscr{B }}
\nc{\ew}{\mathscr{E}}
\nc{\bwc}{\bw_{\aw}}
\nc{\GLoid}{{\mathsf{GL}}}
\nc{\exph}{\ul{\exp}(\h)}

\nc{\lambdanabla}{\PT_{\tau}} 

\nc{\QBw}{\lambda} 
\nc{\MBw}{\mu} 
\nc{\MBwtb}{\PT_{\tau,\bw}} 

\nc{\QPw}{\lambda} 
\nc{\MPw}{\mathscr{P}_{{\veps},\mathscr{B}}} 

\nc{\noQPw}{\lambda_{\veps,\bw}} 
\nc{\noQPwO}{\lambda_{\veps,\bw}} 
\nc{\noMPwO}{\mathscr{P}} 

\nc{\QPwD}{\mathscr{K}} 
\nc{\MPwD}{\mathscr{P}_{\bw}} 

\nc{\tfV}{\Vect_{\hbar}}

\nc{\EK}[2]{\Psi_{#1}^{#2}}
\nc{\FF}[1]{\mathsf{f}_{#1}} 
\nc{\FEK}[1]{\mathsf{F}^{#1}} 

\nc{\FFp}[1]{\mathsf{f}^{\pp}_{#1}} 

\nc{\Mns}[1]{\mathsf{Mns}(#1)}
\nc{\supp}{\mathsf{supp}}
\nc{\zsupp}{\z\supp}
\nc{\DCPS}[1]{G_{#1}}
\nc{\DCPA}[2]{\Upsilon_{#1#2}}
\nc{\DCPAC}[3]{\Upsilon^{#1}_{#2#3}}

\nc{\CUOhinfDI}[2]{\hext{(U\g}^{\Oinf})^{#2}_{\;#1}} 

\nc{\redasso}[2]{\mathsf{a}^{#1}_{#2}}
\nc{\cCox}[1]{\mathscr{C}_{#1}} 
\nc{\CoxS}[3]{S_{#3 #2}^{#1}} 
\nc {\Fi}{F_{\{i\}}}
\nc{\Ki}{\K[i]}
\nc{\Kip}{\K'[i]}
\nc{\Kj}{\K[j]}
\nc{\trunc}[3]{{#1}^{#2}_{#3}}
\nc{\mCox}[1]{\mathbf{H}_{#1}} 
\nc{\nCox}[1]{\mathbf{v}_{#1}} 
\nc{\cc}[1]{\mathsf{conn}(#1)}

\nc{\hRes}{\operatorname{Res}^\hbar}
\nc{\OCox}[2]{\mathscr{O}_{#1}^{#2}}

\nc{\sffh}{{\mathsf f}_\hbar}

\nc{\BBm}{\mathcal{B}_B^{\ulm}}
\nc{\brac}[1]{\mathsf{br}_{#1}}
\nc{\topT}[1]{T_{#1}} 
\rnc{\SS}{\mathfrak{S}}

\nc{\gau}[1]{{\mathsf g}_{#1}}

\nc{\ssmod}{\operatorname{mod}^{\scsop{ss}}}

\nc{\Wh}{\operatorname{Mod}_\h^\hbar}
\nc{\Whp}{\W_\hbar}

\nc {\Pg}{\mathcal{P}_\g}
\nc {\Phg}{\mathcal{P}_\g}
\nc {\Fint}{F}
\nc {\nEK}{Etingof--Kazhdan }
\nc {\qWg}{quantum Weyl group }
\nc {\wrt}{with respect to }
\nc {\Ohint}{\Oh\int}
\nc {\Onhinfintg}{{\mathcal O}^{\hbar,\int}_\g}
\nc {\fml}{[\negthinspace[\hbar]\negthinspace]}
\nc {\Ug}{U\g}
\nc {\Pb}{\PT_{\bw}}
\nc {\Ptb}{\PT_{\tau,\bw}}
\nc {\Pte}{\PT_{\wt{\veps},\bw}}
\nc {\fd}{finite--dimensional }
\nc {\ttg}{\mathcal{T}_\g}
\nc {\res}{\operatorname{Res}}
\nc {\resh}{\operatorname{Res}_\hbar}
\nc {\vepsh}{\veps_\hbar}

\nc {\DCP}{De Concini--Procesi }

\renewcommand {\Im}{\operatorname{Im}}
\newcommand {\rhs}{right--hand side }

\newcommand {\DY}{Drinfeld--Yetter }

\nc{\WJ}{W_{\bfJ}}
\nc{\PB}{\mathcal{PB}}
\nc{\PBJ}{\mathcal{PB}_{\bfJ}}
\nc{\qPOJ}{\O_{\infty,\Uhg}^{\scsop{\bfJ-int}}}
\nc{\qPO}[1]{\O_{\infty,\Uhg}^{\scsop{#1int}}}
\nc{\FFJ}[2]{\mathsf{f}_{#1}^{\scsop{#2int}}}

\nc{\POJ}{\O_{\infty,\g}^{\hbar, \scsop{\bfJ-int}}}
\nc{\PO}[1]{\O_{\infty,\g}^{\scsop{#1int}}}

\nc {\WJPoid}[1]{W_\bfJ\ltimes\Poid{#1}}

\nc{\PJs}{\mathscr{P}_{{\tau}, \bw}} 
\nc{\QPBwJ}{\QPw} 
\nc{\QBwJ}{\QPw^{\scsop{\bfJ--int}}}

\nc{\grpd}{}

\newcommand {\tm}[1]{#1}
\newcommand {\Pwab}{\Pw^{\scriptscriptstyle{\operatorname{ab}}}}

\AtBeginDocument{%
	\def\MR#1{}
}
\author[A.~Appel]{Andrea Appel}
\address{Dipartimento di Scienze Matematiche, 
Fisiche e Informatiche, Universit\`a di Parma, 
and INFN Gruppo Collegato di Parma, 43100 Parma PR, Italy}
\email{\href{mailto:andrea.appel@unipr.it}{andrea.appel@unipr.it}}
\author[V.~Toledano Laredo]{Valerio Toledano Laredo}
\address{Department of Mathematics, Northeastern University, 360 Huntington Avenue, 02115 Boston MA, USA.}
\email{\href{mailto:v.toledanolaredo@northeastern.edu}{v.toledanolaredo@northeastern.edu}}
\title{Pure braid group actions on category $\O$ modules}
\thanks{The first author is partially supported by the Program {\em FIL} 2020 of the University
of Parma and co-sponsored by the Fondazione Cariparma. The second author is partially
supported by the NSF grant DMS--1802412.}
\dedicatory{To Corrado De Concini}

\begin{document}
\begin{abstract}
Let $\g$ be a symmetrisable \KMA and $\Uhg$ its quantised enveloping algebra. Answering
a question of P. Etingof, we prove that the quantum Weyl group operators of $\Uhg$ give rise
to a canonical action of the pure braid group of $\g$ on any category $\O$ (not necessarily integrable)
$\Uhg$--module $\V$. By relying on our recent results \cite{appel-toledano-15}, we show that
this action describes the monodromy of the rational Casimir connection on the $\g$--module
$V$ corresponding to $\V$. 
We also extend these results to yield equivalent 
representations of parabolic pure braid groups on parabolic category $\O$ for \mbox{$\Uhg$ and $\g$.}
\end{abstract}
\Omit{ArXiv abstract
		Let g be a symmetrisable Kac-Moody algebra and U_h(g) its quantised enveloping algebra. Answering
		a question of P. Etingof, we prove that the quantum Weyl group operators of U_h(g) give rise
		to a canonical action of the pure braid group of g on any category O (not necessarily integrable)
		U_h(g)-module V. By relying on our recent results in arXiv:1512.03041, we show that
		this action describes the monodromy of the rational Casimir connection on the g-module
		corresponding to V under the Etingof-Kazhdan equivalence of category O for g and U_h(g).
}
\Omit{2nd ArXiv abstract
		Let g be a symmetrisable Kac-Moody algebra and U_h(g) its quantised enveloping algebra. Answering
		a question of P. Etingof, we prove that the quantum Weyl group operators of U_h(g) give rise
		to a canonical action of the pure braid group of g on any category O (not necessarily integrable)
		U_h(g)-module V. By relying on our recent results in arXiv:1512.03041, we show that
		this action describes the monodromy of the rational Casimir connection on the g-module
		corresponding to V under the Etingof-Kazhdan equivalence of category O for g and U_h(g).
		We also extend these results to yield equivalent quantum Weyl group and monodromic
		representations of parabolic pure braid groups on parabolic category O for U_h(g) and g.
Detailed comment for 2nd arXiv revision
The following changes were made:
Sect 1. The exposition in paragraphs 1.5 and 1.6 was improved. Paragraph 1.12 is new
Sect 2. The definition of category O_infty in 2.2 and 2.3 was corrected so that it is stable under restriction to Levi subalgebras
Sect 3. A remark was added in Sect. 3.1 explaining the relation between the Drinfeld algebra defined in this paper and the algebra originally introduced by Drinfeld
Sect 4. Exposition in 4.7 improved
Sect 5. Exposition in 5.3, 5.4 and 5.5 improved
Sect 8. A footnote was added in 8.1 explaining how our results complete Etingof-Kazhdan's answer to a question raised by Drinfeld.
Sect 9. is a new section which extends the results of the paper to the case of parabolic braid groups
Summary:
1) Corrected the def. of category O_infty in 2.2-2.3 so that it is stable under restriction to Levi subalgebras
2) Added a rk in 3.1 and a footnote in 8.1 explaining the relation between the Drinfeld algebra defined in this
paper and the one originally introduced by Drinfeld 3) A new section 9 was added, which extends our results
to the case of parabolic braid groups 4) Exposition improved throughout
}
	\maketitle\thispagestyle{empty}
	\setcounter{tocdepth}{1}
	\tableofcontents

\section{Introduction}

\subsection{}

Let $\g$ be a symmetrisable \KM algebra, $\Uhg$ its quantized enveloping algebra and
$W$ their Weyl group. We denote by $\O$ the category of deformation highest weight
modules of $\g$, by $\Oint\subset\O$ the full subcategory of integrable ones, and by
$\Ohint\subset\Oh$ the corresponding categories for $\Uhg$. In \cite{appel-toledano-15},
we constructed an equivalence $\Oint\to\Ohint$ which intertwines the
monodromy of the rational Casimir
connection of $\g$ and the quantum Weyl group action of the braid group $\Bw$ of $\g$,
respectively, thus extending the equivalence obtained in \cite{toledano-02,toledano-08,toledano-16}
when $\g$ is finite--dimensional. P. Etingof asked whether this equivalence extends to suitable categories of
modules which are not necessarily integrable, while remaining equivariant under the pure
braid group $\Pw$ of $\g$.

The goal of the present paper is to answer this question in the affermative. Specifically,
we prove that the \qWg action of $\Pw$ on category $\Ohint$ modules can be extended
to all category $\Oh$ modules. We then show that this action is equivalent to the restriction
to $\Pw$ of the equivariant monodromy of the Casimir connection, which is defined on any category
$\O$ module. Our results hold more generally for the category $\O_\infty$ of modules which are
locally finite under the action of the Borel subalgebra, though for simplicity
we restrict to category $\O$ in the Introduction.

\subsection{}

We turn now to a more detailed description of our results. Endow $\O$ with the associativity
and commutativity constraints arising from the KZ equations \cite{drinfeld-90}. In \cite
{etingof-kazhdan-96, etingof-kazhdan-98,
etingof-kazhdan-08}, \nEK constructed a braided tensor equivalence $\FEK{}:\O\to\Oh$
which is Tannakian, that is endowed with a natural isomorphism $\alpha$ fitting in
the diagram 
\begin{equation}\label{eq:intro-diag-1}
	\begin{tikzcd}
		\O \arrow[rr,"\FEK{}"]\arrow[dr,"\FF{}"',""{name=U}]
		&
		&\Oh \arrow[dl,"\FF{\hbar}"]
		\arrow[Rightarrow, sloped, to=U, "\alpha"']\\
		&\Vecth&
	\end{tikzcd}
\end{equation}
where $\Vecth$ is the category of topologically free modules over $\IC\fml$, 
$\sff,\sff_\hbar$ are the forgetful functors, and $\sff$ is endowed with an appropriate
tensor structure. The pair $(\FEK{},\alpha)$ gives rise to an isomorphism $\Psi_\alpha:\End
(\sffh)\to\End(\sff)$ via the composition
\[\End(\sffh)\longrightarrow{} \End(\sffh\circ\FEK{}) \xrightarrow{} \End(\sff)\]
where the first isomorphism is induced by $\FEK{}$, and the second is given by $\Ad(\alpha)$.
Note that $\alpha$ is only unique up to an automorphism $\gamma$ of $\sff$, and that $\Psi
_{\gamma\circ\alpha}=\Ad(\gamma)\circ\Psi_\alpha$.

\subsection{}

Building on our earlier work \cite{appel-toledano-18,appel-toledano-19b,appel-toledano-19a},
we constructed in \cite{appel-toledano-15} an automorphism $\gamma\in\Aut(\sff)$ such
that 
$\Psi_{\gamma\circ\alpha}$ is equivariant \wrt the action of the braid
group $\Bw$ on integrable category $\O$ modules. Specifically, the \nEK functor $\FEK{}$
restricts to an equivalence $\Oint\to\Ohint$
and therefore leads to an isomorphism $\EK{\alpha'}{\sint}:\End(\sffh^{\sint})\to\End(\sff^{\sint})$
for any $\alpha':\FF{\hbar}\circ\FEK{}\Rightarrow\FF{}$. Regard the quantum Weyl group
action of $\Bw$ on objects in $\Ohint$ as a morphism $\QBw:\Bw\to\End(\sff_\hbar^{\sint})$,
and the monodromy of the Casimir connection as a morphism $\MBw:\Bw\to\End(\sff^{\sint})$.
Then, $\gamma$ may be chosen so that the following is a commutative triangle \cite{appel-toledano-15}
\begin{equation}\label{eq:intro-diag-2}
	\begin{tikzcd}
		& \Bw \arrow[dl, " \QBw"'] \arrow[dr, "\MBw"]  & \\
%
		\End(\sff^{\sint}_\hbar)\arrow[rr,"\EK{\gamma\circ\alpha}{\sint}"']
		&
		&\End(\sff^{\sint})
	\end{tikzcd}
\end{equation}
As a consequence, the monodromy of the Casimir connection on a module $V\in\Oint$
is equivalent to the quantum Weyl group action of $\Bw$ on $\FEK{}(V)$.

\subsection{}

P. Etingof asked us whether such an equivalence holds for a larger class of not necessarily
integrable modules, provided $\Bw$ is replaced by the pure braid group $\Pw$. The choice
of the latter is suggested by the fact that $\Bw$ does not act on all category $\O$ modules
for either $\g$ or $\Uhg$, while $\Pw$ does on category $\O$ $\g$--modules via the monodromy
of the Casimir connection.

To the best of our knowledge, no action of $\Pw$ on category $\Oh$ modules has
been previously constructed. The main result of the present paper is to construct such an
action, and then show the commutativity of the resulting diagram
\begin{equation}\label{eq:intro-diag-3}
	\begin{tikzcd}[row sep=tiny]
		& \Pw \arrow[dl, "\lambda"'] \arrow[dr, "\MBw"] & \\
		\End(\sff_\hbar)\arrow[rr,"\EK{\gamma\circ\alpha}{}"']
		&
		&\End(\sff)
	\end{tikzcd}
\end{equation}

\subsection{}\label{ss:intro sign} 

To state our results in more detail, recall first that the abelianisation $\Pwab=\Pw/[\Pw,\Pw]$
of the pure braid group is isomorphic to the free abelian group with a generator $p_\alpha$
for each positive real root $\alpha$ \cite{tits-66,digne-15}. Set $\iota=\sqrt{-1}$, and define
the sign character to be the morphism
\begin{equation}\label{eq:sign intro}
\veps_\hbar:\Pwab\to\Aut(\sff_\hbar^{\sint})
\qquad\qquad
p_\alpha\to\exp(\pi\iota h_\alpha)
\end{equation}
where $\exp(\pi\iota h_\alpha)$ acts as multiplication by $\exp(\pi\iota\nu(h_\alpha))$ on
the $\nu$--weight space of an integrable category $\O_\hbar$ module. The morphism
$\veps_\hbar$ arises as the reduction mod $\hbar$ of the \qWg action of $\Pw$ on
category $\Ohint$.

As a subgroup of $\Bw$, $\Pw$ is generated by the elements $S^2_{w,i}=S_w
S_i^2S_w^{-1}$, where $S_i$ is a generator of $\Bw$, $w\in W$ is such that $w\alpha_
i$ is a positive root, and $S_w\in\Bw$ is the canonical lift of $w$ 
\cite{digne-gomi-01}. Moreover, the \qWg action of $S_{w,i}$ on 
a module $\V\in\Ohint$ is \mbox{given by}
\begin{equation}\label{eq:q Casimir}
\lambda(S^2_{w,i})=
\tm{\exp(\pi\iota h_{w\alpha_i})q^{\ku{w,i}}=\veps_\hbar(S_{w,i})q^{\ku{w,i}}}
\end{equation}
where the second factor is the truncated quantum Casimir operator for the copy
of $U_\hbar\sl{2}\subset\Uhg$ corresponding to the pair $(w,i)$ \cite{lusztig-93},
and $q=\exp(\hbar/2)$.

\subsection{}\label{ss:extend to Oh}

To extend this action to an arbitrary category $\Oh$ module, we lift the sign
character $\veps_\hbar$ to a morphism 
\[\Pwab\to\Aut(\sff_\hbar)
\qquad\qquad 
p_\alpha\to\exp(\pi\iota h_{\alpha})\]
which we denote by the same symbol. We then prove that the quantum
Casimirs $\tm{q^{\ku{w,i}}\in\Uhg}$  
give rise to a morphism $\QPwD:\Pw\to
(\Uhg)^\h$. It follows that
\begin{equation}\label{eq:lambda intro}
\QPw:\Pw\to\Aut(\sff_\hbar)
\qquad\qquad
S_{w,i}^2 \to \exp(\pi\iota h_{w\alpha_i})q^{\ku{w,i}}
\end{equation}
is an extension of the \qWg action of $\Pw$ to all category $\Oh$ modules.

\subsection{} 

The fact that $\QPwD$ is a morphism would follow at once if $\End(\sff_\hbar)$
acted faithfully on $\sff_\hbar^{\sint}$. This, however, is clearly false: if $\varphi$
is any function on $\h^*$ which vanishes on integral weights, then $\varphi\in
\End(\sff_\hbar)$, but $\varphi$ maps to zero in $\End(\sff_\hbar^{\sint})$.
To remedy this, we rely on the fact that $\Uhg$ acts faithfully on $\sff\int_\hbar$,
whose proof is due to Etingof. This implies that any $\lambda(p)\in\End(\sff_
\hbar^{\sint})$, $p\in\Pw$, arises from the action of a unique element of $\Uhg$,
thereby yielding the required action of $\Pw$ on $\End(\sff_\hbar)$.\footnote{Note
that this bypasses having to explicitly check that the quantum Casimirs satisfy
the relations of the generators $S_{w,i}$ given in \cite[Cor. 6]{digne-gomi-01}.}

A similar argument works for the quantum group $\Uqg$, where $q$ is either
an indeterminate, or not a root of unity. In that case, the quantum Casimirs
$\tm{q^{\ku{w,i}}}$ 
do not lie in $\Uqg$, but in a variant $\DrAq$ of an algebra originally introduced
by Drinfeld \cite[Sect. 8]{drinfeld-problem-92}, which consists of formal, infinite
series of the form $\sum c_{X}X$, where $X$ runs over a weight basis of $\Uqnp
{}$ and $c_X\in\Uqbm{}$. Etingof's faithfulness result also applies to $\DrAq$,
and yields an action of $\Pw$ on any category $\O$ module for $\Uqg$.

\subsection{} 

Let now $\sfY$ be the complexification of the Tits cone of $\g$, $\sfX\subset\sfY$ 
its set of regular points, and $x_0\in\sfX$ a basepoint. By a theorem of van der Lek
\cite{van-der-lek-83}, which generalises Brieskorn's \cite{brieskorn-71}, the pure
and full braid groups may be realised as
\[\Pw\cong\Poid(\sfX;x_0)\aand\Bw\cong\Poid(\sfX/W;[x_0])\]

The Casimir connection is the $\Ug$--valued formal meromorphic connection
on $\sfX$ with logarithmic singularities on the root hyperplanes given by 
\begin{equation}\label{eq:nabla intro}
\nablak=
d-\nablah\sum_{\alpha\succ 0}
\frac{d\alpha}{\alpha}\cdot\Ku{\alpha}{+}
\end{equation}
where $\Ku{\alpha}{+}=\sum_{i=1}^{\rsm{\alpha}}e_{-\alpha}^{(i)}e_\alpha^{(i)}$
is the normally ordered truncated Casimir operator corresponding to the positive
root $\alpha$, and $\sfh=\hbar/2\pi\iota$ \cite{millson-toledano-05,
toledano-02,procesi-96,felder-markov-tarasov-varchenko}.
The sum \eqref{eq:nabla intro} over $\alpha$ is locally finite on any (not necessarily
integrable) category $\O$ module $V$, and gives rise to a well--defined flat connection
on the holomorphically trivial vector bundle $\IV$ on $\sfY$ with fibre $V$. Its monodromy
therefore gives rise to a morphism
\begin{equation}\label{eq:PT}
\PT:\Poid(\sfX;x_0) \to \End(\sff) 
\end{equation} 

\subsection{} \label{ss:intro equiv} 

The normal ordering in \eqref{eq:nabla intro} breaks the equivariance of 
$\nablak$ \wrt the action of $W$ on $\sfX$ and the subalgebra of $\h$--invariants
$\Ug^\h\subset\Ug$, which contains the Casimirs $\Ku{\alpha}{+}$.

Nevertheless, it is possible to modify the monodromy of $\nablak$ so that it gives
rise to a representation of the braid group $\Bw$ on integrable category $\O$ modules
\cite[Sect. 4]{appel-toledano-15} (see also Section~\ref{s:Casimir}). This relies on
the equivalence of groupoids
\begin{equation}\label{eq:gp to gpd}
\cE_{x_0}:\Poid(\sfX/W;[x_0])
\to W\ltimes \Poid(\sfX;Wx_0)
\end{equation}
where the \rhs is the semi--direct product of $W$ with the fundamental groupoid of
$\sfX$ based at the orbit $Wx_0$, and $\cE_{x_0}$ is given by the unique lifting of
loops through $x_0$, and proceeds as follows.
\begin{itemize}[leftmargin=2em]\itemsep0.1cm
\item Extend the monodromy of $\nablak$ to a morphism
\begin{equation}\label{eq:PT intro}
\PT:\Poid(\sfX;Wx_0)\to\End(\sff)
\end{equation}
\item Replace the target of $\PT$ by a subalgebra
$\ttg\subset\End(\sff)$ which, unlike $\End(\sff)$, is acted upon by $W$.
$\ttg$ is the image of the holonomy algebra of the root arrangement of
$\g$, and is a completion of the subalgebra of $\Ug^\h\fml$ generated
by the Casimirs $\hbar\Ku{\alpha}{+}$ and the Cartan subalgebra $\hbar\h$.
\item The lack of equivariance of $\nablak$ can then be measured by a
1--cocycle
\[\aw:W\to\Hom_{\grpd}(\Poid{(\sfX;Wx_0)},\ttg)\]
defined by
$\aw_w(\gamma)=\PT(\gamma)^{-1}\cdot w^{-1}\PT(w\gamma)$.
\item We prove that $\aw$ is abelian \ie takes values in 
\[\sfM=\Hom_{\grpd}(\Poid{(\sfX;Wx_0)},\exp(\hbar\h))\]
and that it is the coboundary of an essentially unique cochain
$\bw\in\sfM$ \ie that $\aw_w=\bw\cdot(w^{-1}\bw)^{-1}$ for any
$w\in W$.
\item As a consequence, $\PT$ can be modified to a $W$--equivariant
morphism 
\[\Pb:\Poid(\sfX;Wx_0)\to\ttg
\qquad\qquad
\Pb(\gamma)=\PT(\gamma)\cdot\bw(\gamma)\]
\item Composing $\Pb$ with the morphism $\cE_{x_0}$ 
\eqref{eq:gp to gpd} then yields an action of $\Bw$ on any $W\ltimes\ttg
$--module.
\item It is well--known that $W$ does not act on an integrable module
$V$, but that the triple exponentials
\begin{equation}\label{eq:tau}
\tau_i=
\exp(e_i)\cdot\exp(-f_i)\cdot\exp(e_i)
\end{equation}
are well--defined on $V$, permute its weight spaces according to the
$W$--action, and give rise to a morphism $\tau:\Bw\to\Aut(\sff\int)$.
\item Finally, lifting $\cE_{x_0}$ to
$\wt{\cE}_{x_0}:\Poid(\sfX/W;[x_0])
\to \Bw\ltimes \Poid(\sfX;Wx_0)$,
and composing with $\tau\otimes\Pb$ yields a morphism
\begin{equation}\label{eq:intro equi mon}
\Ptb:\Bw\to\Aut(\sff\int)
\qquad\qquad 
\gamma \to \tau(\gamma)\cdot\PT(\gamma)\cdot\bw(\gamma)
\end{equation}
which we term the {\it equivariant monodromy} of $\nablak$.
\end{itemize}

\subsection{} 

By \cite{appel-toledano-15}, the equivariant monodromy of $\nablak$ on an
integrable module $V\in\Oint$ is canonically equivalent to the \qWg action of
$\Bw$ on the \nEK quantisation $\FEK{}(V)\in\Ohint$ \ie the diagram \eqref
{eq:intro-diag-2} is commutative for $\MBw=\Ptb$. This can be used to give
a monodromic description of the action $\lambda$ \eqref{eq:lambda intro}
of $\Pw$ on category $\Oh$ modules as follows. 

The restriction of the triple exponential map $\tau$ \eqref{eq:tau}
to $\Pw$ is the sign character 
\[\veps:\Pwab\to\Aut(\sff^{\sint})
\qquad\qquad
p_\alpha\to\exp(\pi\iota h_\alpha)\]
Lifting it to $\veps:\Pwab\to\Aut(\sff)$ as in \ref{ss:extend to Oh} therefore
lifts the equivariant monodromy action of $\Pw$ to
\begin{equation}\label{eq:Pte}
\MPw:\Pw\to\Aut(\sff)
\qquad\qquad
\gamma \to{\veps}(\gamma)\cdot\PT(\gamma)\cdot\bw(\gamma)
\end{equation}
\ie extends the restriction of $\Ptb$ to $\Pw$ to any category $\O$ module.

\nc {\bbeta}{\gamma\circ\alpha}
To relate $\MPw$ to $\lambda$, denote the restriction morphisms by
\[\res:\End(\sff)\to\End(\sff\int)\aand
\resh:\End(\sff)\to\End(\sff\int)\]
The commutativity of \eqref{eq:intro-diag-2} implies that, for any $p\in\Pw$
\[
\res\circ\EK{\bbeta}{}\circ\QPw(p)=
\EK{\bbeta}{\sint}\circ\resh\circ\QPw(p)=
\res\circ\MPw(p)\]
and therefore that $\res\circ\EK{\bbeta}{}\circ\QPwD(p)=\res\circ(\PT(p)\bw(p))$,
since ${\veps}=\EK{\bbeta}{}({\veps}_\hbar)$. In turn, this
implies that $\EK{\bbeta}{}\circ\QPwD(p)=\PT(p)\bw(p)$, so that $\EK{\bbeta}{}$
intertwines $\QPw$ and $\MPw$, since $\EK{\bbeta}{}$ maps the Drinfeld algebra
$\DrAh{}\supset\Uhg$ to its classical analogue $\DrA{}$, the latter acts faithfully on $\sff$,
and the algebra $\ttg\ni\PT(p),\bw(p)$ is contained in $\DrA{}$.

\subsection{}\label{ss:no-qWg} 

The above can also be used to give a description of the (non--equivariant)
monodromy $\PT:\Pw\to\Aut(\sff)$ of the Casimir connection $\nablak$ \eqref
{eq:PT} in terms of quantum Weyl group operators as follows.

We prove that the restriction to $\Pw$ of the cochain $\bw$ is the map $\Pwab
\to\exp(\hbar\h)$ given by $\bw(p_\alpha)=\exp(\hbar\tinv{\alpha}/2)$, 
where $\tinv{\alpha}\in\h$ corresponds to $\alpha$ via the isomorphism $\h^
*\to\h$ induced by the chosen inner product on $\g$. Define the morphism
\[\noQPw:\Pw\to\Aut(\sffh) 
\qquad\qquad
p \to 
\vepsh(p)^{-1}\cdot \lambda(p)\cdot \bw(p)^{-1}=
\QPwD(p)\cdot\bw(p)^{-1}\]
We refer to $\noQPw$ as the {\it normally ordered} \qWg action of $\Pw$
on category $\Oh$ modules. The terminology is motivated by the fact that,
while $\lambda(S_i^2)=\exp(\pi\iota h_i)\cdot q^{\tm{\ku{\hbar,i}}}$ 
by \eqref{eq:lambda intro}, $\noQPw(S_i^2)=q^{2\tm{\ku{\hbar,i}^+}}$,
where the latter is a normally ordered version of the quantum Casimir.
The commutativity of \eqref{eq:intro-diag-2} then implies
that $\noQPwO$ computes the monodromy of $\nablak$, that is that
$\EK{\bbeta}{}\circ\noQPwO=\noMPwO$. 

\subsection{}\label{ss:intro-parabolic}

The above results can be generalised to the parabolic setting as follows. Let
$\bfJ$ be a subset of nodes of the Dynkin diagram of $\g$, $\g_\bfJ\subseteq
\g$ the corresponding Lie subalgebra, $W_\bfJ\subseteq W$ its Weyl group,
and $\PBJ\subseteq\Bw$ the parabolic pure braid group given by the preimage of
$W_\bfJ$.

We construct a \qWg action of $\PBJ$ on any category $\O_\hbar$ module whose
restriction to $\Uhg_\bfJ$ is integrable. This action is such that
\begin{itemize}\itemsep0.25cm
	\item its restriction to the braid group $\Br{W_\bfJ}$ is the \qWg 
	action of $\Br{W_\bfJ}$ on integrable $\Uhg_\bfJ$--modules
	\item its restriction to the pure braid group $\Pw$ coincides with the \qWg
	action \eqref{eq:lambda intro} on category $\O_\hbar$ modules
\end{itemize}
We also define a normally ordered version of this \qWg action, in analogy with \ref
{ss:no-qWg}.

We then construct a monodromy action of $\PBJ$ on any category $\O$ module
whose restriction to $\g_\bfJ$ is integrable. We do so by relying on the fact that
$\PBJ$ is isomorphic to $\Poid{(\sfX/W_\bfJ;[x_0])}$, and correcting the equivariance
of the Casimir connection, as outlined in \ref{ss:intro equiv}, but only \wrt $W_\bfJ$.
The resulting $W_\bfJ$--equivariant monodromy action is such that
\begin{itemize}\itemsep0.25cm
\item its restriction to $\Br{W_\bfJ}$ is the equivariant monodromy action
of $\Br{W_\bfJ}$ on integrable category $\O$ $\g_\bfJ$--modules
\item its restriction to $\Pw$ coincides with the monodromy action \eqref{eq:PT intro}
on category $\O$ modules (up to a simple correction on $\P_{W_\bfJ}$).
\end{itemize}

Finally, we show that the above \qWg and monodromic actions of $\PBJ$ are equivalent
by relying on the fact that $\PBJ$ is generated by $\Br{W_\bfJ}$ and $\Pw$, and using
the equivalences \eqref{eq:intro-diag-2} for $\Br{W_\bfJ}$ and \eqref{eq:intro-diag-3} for $\Pw$.

\subsection{Outline of the paper} 

In Section \ref{s:km}, we review the definition of quantum Weyl group operators. In Section
\ref{s:faithful}, we introduce the Drinfeld algebra and prove that it acts faithfully on $\Ohint$.
In Section \ref{s:pure-qWg}, we construct the \qWg action of $\Pw$ on category $\O$. Section
\ref{s:Casimir} reviews the definition of the Casimir connection, and the equivariant extension
of its monodromy to a representation of the braid group $\Bw$. Section \ref{s:mono-recap-1}
reviews the definition of a braided Coxeter category, and Section \ref{s:mono-recap-2} the
main result of \cite{appel-toledano-15}. In Section \ref{s:main}, we prove the stated equivalence.
We also point out that it continues to hold if $\FEK{}$ is replaced by the \nEK equivalence
$\FEK{\sfPhi}$ corresponding to an arbitrary Lie associator $\sfPhi$ rather than the one arising
from the KZ equations. Finally, in Section \ref{s:parabolic}, we generalise these results to
parabolic pure braid groups.

\subsection{Acknowledgments} 
We heartily thank Pavel Etingof for persistently asking us about extending the main
result of \cite{toledano-08,toledano-16,appel-toledano-15} to pure braid groups, and
for his detailed comments on, and very careful reading of, an earlier version of this paper.
\newpage

\section{\KM algebras and quantum groups}\label{s:km}

\subsection{Symmetrisable Kac--Moody algebras \cite{kac-90}}\label{ss:sym-km-summary}

Let $\bfI$ be a finite set and $\GCM{A}=(a_{ij})_{i,j\in\bfI}$ a generalised Cartan matrix, \ie
$a_{ii}=2$, $a_{ij}\in\IZ_{\leqslant 0}$, $i\neq j$, and $a_{ij}=0$ implies $a_{ji}=0$.
Let $(\h,\Pi,\Pi^\vee)$ be a realization of $\GCM{A}$, \ie 
\begin{itemize}\itemsep0.25cm
	\item $\h$ is a finite--dimensional complex vector space\footnote{Note that, unlike \cite{kac-90},
	we do not require $\h$ to have minimal dimension
	$2|\bfI|-\operatorname{rank}(A)$.} 
	\item $\Pi=\{\alpha_i\}_{i\in\bfI}$ is a linearly independent subset of $\h^*$
	\item $\Pi^\vee=\{\crt{i}\}_{i\in\bfI}$ is a linearly independent subset of $\h$
	\item $\alpha_i(\crt{j})=a_{ji}$ for any $i,j\in\bfI$ 
\end{itemize}

The Kac--Moody algebra corresponding to $\GCM{A}$ and the realisation $(\h,\Pi,\Pi^\vee)$
is the Lie algebra $\g$ generated by $\h$ and elements $\{e_i, f_i\}_{i\in\bfI}$, with relations $[\h,
\h]=0$ and
\begin{gather}
	[h,e_i]=\alpha_i(h)e_i
	\qquad\qquad
	[h,f_i]=-\alpha_i(h)f_i
	\qquad\qquad
	[e_i,f_j]=\drc{ij}\crt{i}
\end{gather}
and, for any $i\neq j$,
\[
\mathsf{ad}(e_i)^{1-a_{ij}}(e_j)=0=\mathsf{ad}(f_i)^{1-a_{ij}}(f_j)
\]
Let $\n^\pm\subset\g$ be the Lie subalgebras generated
by $\{e_i\}_{i\in\bfI}$ and $\{f_i\}_{i\in\bfI}$, respectively.

Assume that $\GCM{A}$ is symmetrisable, and fix an invertible diagonal matrix $\GCM{D}=
\Diag(\symd{i})_{i\in\bfI}$ with coprime entries $\symd{i}\in\IZ_{>0}$ such that $\GCM{DA}$
is symmetric. Then, there is a symmetric, non--degenerate bilinear form $\iip{\cdot}{\cdot}$
on $\h$ such that $\iip{\crt{i}}{-}=\symdi{i}\alpha_i$ (see, \eg \cite[Prop.~11.4]{appel-toledano-19a}).
The corresponding identification $\nu:\h\to\h^*$ intertwines the actions of $W$, 
satisfies $\nu(\crt{i})=\symdi{i}\alpha_i$ and therefore restricts to an isomorphism $\h'\stackrel
{\sim}{\longrightarrow}\sfQ\ten_{\IZ}\IC$, where $\h'$ is the span of $\{h_i\}_{i\in\bfI}$ and
${\mathsf{Q}}=\bigoplus_{i\in\bfI}\IZ{\alpha}_i\subseteq{\h}^*$ is the root lattice. Note that
$\iip{\crt{i}}{\crt{i}}=2\symdi{i}$, while the induced form  on $\h^*$ satisfies $\iip{\rt{i}}{\rt{i}}
=2\symd{i}\in 2\IZ_{>0}$. 

By \cite[Thm.~2.2]{kac-90}, $\iip{\cdot}{\cdot}$ uniquely extends to a non--degenerate,
invariant symmetric bilinear form on $\g$, which satisfies $\iip{e_i}{f_j}=\drc{ij}\symdi{i}$
and $[x,y]=\iip{x}{y}\cdot\tinv{\alpha}$ for any $x\in\g_{\alpha}$, $y\in\g_{-\alpha}$, where
$\tinv{\alpha}=\nu^{-1}(\alpha)$.

\subsection{Category $\Oinf$ representations}\label{ss:cat-O-diag}

If $V$ is an $\h$--module and $\lambda\in\h^*$, we denote the corresponding
weight space of $V$ by
\[V[\lambda]=\{v\in V|\,h\,v=\lambda(h)v,\,h\in\h\}\]
and set $P(V)=\{\lambda\in\h^*|\,V{[\lambda]}\neq0\}$.
A $\g$--module $V$ is
\begin{enumerate}[label=(C\arabic*)]\itemsep0.25cm
\item\label{cond:weight-dec} a {\it weight module} if $V=\bigoplus_{\lambda\in\h^*}
V{[\lambda]}$.
\item\label{cond:int} {\it integrable} if it is a weight module, and the elements $\{e_i,
f_i\}_{i\in\bfI}$ act locally nilpotently.
	
This implies that $\lambda(h_i)\in\IZ$ for any $\lambda\in P(V)$ and $i\in\bfI$,
and that $V$ is completely reducible as a (possibly infinite) direct sum of simple
\fd modules over $\sl{2}^{\alpha_i}=\langle e_i, h_i,f_i\rangle\subset\g$.
\item\label{cond:Oinfty} in {\it category $\Oinfg$} if the action of $\bp{}$  is locally finite,
\ie any $v\in V$ is contained in a \fd $\bp{}$--submodule of $V$.
	
This implies in particular that $V$ is the direct sum of its generalised weight spaces
and that, for any $v\in V$, $(U\n^+)_\beta v=0$ for all but finitely many $\beta\in\sfQ
_+$.
\item\label{cond:O} in {\it category $\Og$} if it is a weight module with \fd weight
spaces, such that
\begin{equation}\label{eq:cone condition}
P(V)\subseteq D(\lambda_1)\cup\cdots\cup D(\lambda_m)
\end{equation}
for some $\lambda_1,\ldots,\lambda_m\in\h^*$, where 
$D(\lambda)=\{\mu\in\h^*\;|\;\mu\leqslant\lambda\}$ and
$\mu\leqslant\lambda$ iff $\lambda-\mu\in\sfQ_+=\bigoplus_{i\in\bfI}\bbN\alpha_i$
\end{enumerate}

The categories $\Og\subset\Oinfg$ are symmetric tensor categories. Denoting by $\Ointg\subset
\Og$ and $\Oinfintg\subset\Oinfg$ the full tensor subcategories of integrable representations, we
have the following inclusions
\[\xymatrix@C=.6cm@R=.2cm
{\Og&\subset&\Oinfg\\
\cup&&\cup\\
\Ointg&\subset&\Oinfintg
}\]

\subsection{Deformation category $\Oinf$ representations}\label{ss:cat-O-diag-h}

Similar notions can be defined for $\g$--modules in the category $\Vect_{\hbar}$
of topologically free $\hext{\IC}$--modules. Namely, a $\g$--module $\V\in\Vect_{\hbar}$ 
is called 
\begin{enumerate}[label=(D\arabic*)]\itemsep0.25cm
\item a {\it weight module} if $\V=\bigoplus_{\lambda\in\h^*}\V[\lambda]$,\footnote{Note
that the eigenvalues of the action of $\h$ on $\V$ are required to lie in $\h^*\subsetneq
\h^*\fml$.} where $\bigoplus$ is the direct sum in $\Vect_{\hbar}$, \ie the completion of
the algebraic direct sum in the $\hbar$--adic topology.

\item {\it integrable} if it is a weight module and, for any $i\in\bfI$ and $v\in\V$,
$\lim_{n\to\infty}e_i^n v=0=\lim_{n\to\infty}f^n_iv$.
%
\item\label{cond:Oinftyh} in {\it category $\Ohinfg$} if the action of $\bp{}$ on $\V/\hbar^n\V$
is locally finite for any $n\geq 0$.	
\item\label{cond:O-h} in {\it category $\Ohg$} if it is a weight representation with finite--rank
weight spaces, and such that $P(\V)$ satisfies \eqref{eq:cone condition}.
\end{enumerate}

It is easy to see that $\V$ is a weight (resp. integrable) module in $\Vect_{\hbar}$
if and only if $\V/\hbar^n\V$ is a weight (resp. integrable) module in $\Vect$
for any $n\geq 0$.
We denote by $\Ohintg\subset\Ohg$ and $\Ohinfintg\subset\Ohinfg$ 
the full tensor subcategories of integrable representations.

\subsection{Braid group action}\label{ss:km-triple}

Let $W$ be the Weyl group of $\g$, and $\{s_i\}_{i\in\bfI}$ its set of simple
reflections. 
The braid group $\Bw$ is the group generated by the elements $\{\topS{i}\}_{i\in\bfI}$,
with relations
\begin{equation}\label{eq:gen-braid}
	\underbrace{\topS{i}\cdot \topS{j}\cdot \topS{i}\;\cdots\;}_{m_{ij}}=
	\underbrace{\topS{j}\cdot \topS{i}\cdot \topS{j}\;\cdots\;}_{m_{ij}}
\end{equation}
for any $i\neq j$, where $m_{ij}$ is the order of $s_is_j$ in $W$.
If $V$ is an integrable $\g$--module in $\Vect$ or $\Vect_{\hbar}$, the operators
\begin{equation}\label{eq:triple}
\texp{i}=\exp(e_i)\cdot\exp(-f_i)\cdot\exp(e_i)\in GL(V)
\end{equation}
are well--defined, and satisfy the braid relations \eqref{eq:gen-braid} \cite{tits-66}.
The corresponding action of $\Br{W}$ on $V$ factors through the {\it Tits extension}
$\wt{W}$, an extension of $W$ by the sign group $\IZ_2^\bfI$ .


\subsection{The quantum group $\Uhg$ \cite{drinfeld-quantum-groups-87,jimbo-85}}\label{ss:DJ-qg}

Let $\hbar$ be a formal variable, set $q=\exp(\hbar/2)$ 
and $q_i= q^{\symd{i}}$, $i\in\bfI$. 
The Drinfeld--Jimbo quantum group of $\g$ is 
the algebra $\Uhg$ over $\hext{\IC}$ topologically generated by 
$\h$ and the elements $\{E_i, F_i\}_{ i\in\bfI}$, subject to the relations $[h,h']=0$, 
\begin{align*}
	[h, E_i]=\alpha_i(h)E_i
	\qquad\qquad
	[h, F_i]=-\alpha_i(h)F_i
	\qquad\qquad
	[E_i, F_j]=\drc{ij}\frac{\qKi{i}{}-\qKi{i}{-}}{q_i-q_i^{-1}}
\end{align*}
for any $h,h'\in\h$, $i,j\in\bfI$, and the $q$--Serre relations
\begin{align}
	\sum_{m=0}^{1-a_{ij}}(-1)^m{1-a_{ij}\brack m}_i X_i^{1-a_{ij}-m}X_jX_i^{m}&=0
\end{align}
for $X=E,F$, $i\neq j\in\bfI$, where $\displaystyle [n]_i=\frac{q_i^n-q_i^{-n}}{q_i-q_i^{-1}}$ 
and, for any $k\leqslant n$,
\begin{align}
	[n]_i!= [n]_i\cdot [n-1]_i\cdots[1]_i
	\aand
	{n\brack k}_i=\frac{[n]_i!}{[k]_i!\cdot [n-k]_i}	
\end{align}



Define weight, integrable, category $\Oinf$ and $\O$ modules for $\Uhg$ in $\Vect_\hbar$
analogously to Section \ref{ss:cat-O-diag-h}, and denote by
\[\OinfintUhg\subset\OinfUhg
\aand
\O_{\Uhg}^{\sint}\subset\OUhg\]
the subcategories of integrable modules.\footnote{note in particular that a representation
$\V$ of $\Uhg$ is in category $\Oinf$ if the action of $\Uhbp{}$ on $\V/\hbar^n\V$ is locally
finite for any $n\geq 0$.} 

\subsection{Quantum Weyl group operators \cite{kirillov-reshetikhin-90,lusztig-90,lusztig-93,sa-94,so-90}}\label{ss:qWg}

For any $\V\in\OinfintUhg$, define the endomorphisms $\{\qWSk{i}\}_{i\in\bfI}$ of $\V$ as
follows.\footnote{The operator $\qWSk{i}$ is the operator $T''_{i,+1}$ defined in \cite
[Sec.~5.2]{lusztig-93}.} For any $v_{\mu}\in \V[\mu]$, set
\begin{equation}\label{eq:qwg-op}
	\qWSk{i}\;v_\mu=
	\sum_{\substack{a,b,c \in \IZ_{\geqslant 0}\\ a-b+c = -\mu(\crt{i})} }
	(-1)^b q_i^{b-ac} E_{i}^{(a)} F_{i}^{(b)} E_{i}^{(c)}\cdot v_{\mu}
\end{equation} 
where $X_i^{(a)}= X_i^a/[a]_i!$.

Then, $\qWSk{i}(\V[\mu]) \subseteq \V[s_i(\mu)]$ and the $\qWSk{i}$ give rise to an action of
the braid group $\Br{W}$ on $\V$, which deforms the action by triple exponentials described
in \ref{ss:km-triple} \cite[Sec.~39.4]{lusztig-93}.

\subsection{Action of $\Bw$ on $\Uhg$ (\cite{lusztig-88}, \cite[Chap. 37--39]{lusztig-93})}\label{ss:qWg-aut}

Consider the algebra automorphisms $\{\qWT{i}\}_{i\in\bfI}$ of $\Uhg$ defined by
\[\qWT{i}(h)=s_i(h)
\qquad\qquad
\qWT{i}(E_i)=-F_iq_i^{\crt{i}}
\qquad\qquad
\qWT{i}(F_i)=-q_i^{-\crt{i}}E_i\]
where $h\in\h$ and, for any $i\neq j\in\bfI$,
\[\qWT{i}(X_j)=\sum_{r=0}^{-a_{ij}}(-1)^rq_i^{{\sigma(X)}r}X_i^{-a_{ij}-r}X_jX_i^r\]
where $X=E,F$ and $\sigma(E)=-1=-\sigma(F)$. 

The automorphisms $\{\qWT{i}\}_{i\in\bfI}$ define an action of the braid group $\Bw$ on $\Uhg$
which we denote by $b(X)$, $b\in\Bw$ and $X\in\Uhg$. Moreover, for any $X\in\Uhg$,
$\V\in\OinfintUhg$, and $v\in \V$, one has $\qWSk{i}(X\cdot v)=\qWT{i}(X)\cdot \qWSk{i}(v)$. 


\section{Faithfulness of category $\O$ integrable modules}\label{s:faithful}

Integrable $\Uhg$-modules are well--known to be faithful, \ie the only element of $\Uhg$
acting trivially on every integrable module is zero \cite[Prop.~3.5.4]{lusztig-93}. To
the best of our knowledge, the analogous result for the more restrictive class of integrable
modules in category $\O$ does not appear in the literature. We present here a proof
due to P. Etingof, which establishes faithfulness for a larger algebra containing $\Uhg$.

\subsection{The Drinfeld algebra $\DrAh$}

\label{ss:quantum-drinfeld-series}\label{ss:quantum-faithful}
For any $\beta\in\Qp$, let $\B_{\beta}=\{X_{\beta,p}\}$ be a basis of $\Uhnp{\beta}$ and
set $\B=\bigsqcup_{\beta\in\Qp}\B_{\beta}$. Set
\begin{equation}\label{eq:drinfeld}
	\DrAh{}=\left
	\{\sum_{X\in\B}c_{X}X\,\colon\, c_{X}\in \Uhbm{}\right\}
	={\prod_{\beta\in\Qp}}
	\Uhbm{}\ten\Uhnp{\beta}\supset \Uhg
\end{equation}

$\DrAh{}$ has an algebra structure which extends that of $\Uhg$. Moreover, the action of
$\Uhg$ on any module $\V\in\OUhg$ extends to one of $\DrAh{}$ since, for any $v\in\V$,
$\Uhnp{\beta}v=0$ for all but finitely many $\beta\in\Qp$. 

\begin{theorem}[Etingof]\label{thm:faithful}
	Category $\O$ integrable $\Uhg$--modules are faithful for $\DrAh{}$.
\end{theorem}

The proof is carried out in Sections \ref{ss:quantum-verma}--\ref{ss:proof-quantum-faithfulness}.

\begin{remark}\label{rk:Qh}
A variant $\DrQh$ of the algebra $\DrAh$ was introduced by Drinfeld in \cite[Sect. 8]{drinfeld-problem-92}
as follows. For any $\beta\in\Qp$, let $I_\beta\subset\Uhg$ be the left ideal generated by
$\Uhg_{\beta'}$ for any $\beta'>\beta$, or equivalently by $\{\Uhnp{\beta'}\}_{\beta'>\beta}$,
and set $\DrQh=\lim_{\beta}\Uhg/I_\beta$. Since $\Uhg/I_\beta\cong\bigoplus_{\beta'\ngtr
\beta}\Uhbm{}\ten\Uhnp{\beta'}$, $\DrQh$ embeds into $\DrAh$ as the subalgebra consisting
of series $\sum_{\beta\in\Qp}X_\beta$, $X_\beta\in\Uhbm{}\ten\Uhnp{\beta}$, where for any
$\beta\in\Qp$, $X_{\beta'}=0$ for all but finitely many $\beta'\ngtr \beta$. 
The algebra $\DrQh$ is less natural than $\DrAh$, however. For instance, if $\emptyset\subsetneq\bfJ
\subsetneq\bfI$ is a proper non--empty subset, $\g_\bfJ\subset\g$ the corresponding subalgebra, and $\Q_{\bfJ,\hbar}$
(resp. $\D_{\bfJ,\hbar}$) the analogue of $\DrQh$ (resp. $\DrAh$) for $\g_\bfJ$, then $\D_
{\bfJ,\hbar}\subset\DrAh$ while $\Q_{\bfJ,\hbar}$ does not map to $\DrQh$.

\end{remark}


\subsection{Verma modules}\label{ss:quantum-verma}
For $\lambda\in{\h^*}$, let $M(\lambda)$ be the Verma module of highest weight $\lambda$
and $v_{\lambda}\in M(\lambda)$ its cyclic vector.
For any $\beta\in\Qp$, let $M(\lambda)_{\beta}\subset M(\lambda)$ be the weight space 
of weight $\lambda-\beta$.
Note that there is a natural identification $M(\lambda)_{\beta}\simeq (\Uhnm{})_{\beta}$.
Recall that the {\em contragredient} Verma module $M^{\vee}(\lambda)$ is the 
pullback through the Chevalley involution of the restricted dual $M^*(\lambda)=\bigoplus_{\beta\in\Qp}M(\lambda)_{\beta}^*$, where
$M(\lambda)_{\beta}^*$ denotes the dual in $\Vecth$.
The contragredient Verma module is  
equipped with a morphism $M(\lambda)\to M^{\vee}(\lambda)$, $v_{\lambda}\mapsto v_{\lambda}^*$.
The Shapovalov form on $M(\lambda)$ is defined by
\begin{align}
	\iip{\cdot}{\cdot}_{\lambda}: M(\lambda)\ten M(\lambda)\to M(\lambda)\ten M^{\vee}(\lambda)\to \hext{\IC}
\end{align}
By construction, it satisfies $\iip{v_{\lambda}}{v_{\lambda}}_{\lambda}=1$, 
$\iip{M(\lambda)_{\beta}}{M(\lambda)_{\beta'}}_{\lambda}=0$ if $\beta\neq\beta'$,
and $\iip{xv}{w}_{\lambda}=-\iip{x}{\omega(x)w}_{\lambda}$ 
for any $x\in\g$, $v,w\in M(\lambda)$. It is well--known that 
$\iip{\cdot}{\cdot}_{\lambda}$ is symmetric and non--degenerate only for 
generic $\lambda\in\h^*$.

For generic $\lambda\in\h^*$, let $\B_{\lambda,\beta}^*=\{X_{\beta,p}^*\}$
be the dual  basis of $\Uhnm{\beta}$ with respect to the Shapovalov form. 
In particular, one has $\iip{X_{\beta,i}^*v_{\lambda}}{\omega(S(X_{\beta,j}))v_{\lambda}}=\drc{ij}$.
Thus, modulo elements of weights lower than $\lambda$,
$X_{\beta,j} X_{\beta,i}^* v_{\lambda}= \drc{ij}v_{\lambda}$.

\begin{proposition}
	Verma modules are faithful for $\DrAh{}$.
\end{proposition}

\begin{proof}
	Assume that $x\in\DrAh{}$ acts trivially on every $M(\lambda)$, and write
	\begin{align}
		x=\sum_{\B}x^-_{\beta, i} x^0_{\beta, i} X_{\beta, i}
	\end{align}
	where $x^0_{\beta, i}\in \hext{U\h}$ and $x^-_{\beta, i}\in \Uhnm{}$. Note that, for any $\lambda\in\h^*$, the action of $x$ on the cyclic vector $v_{\lambda}\in M(\lambda)$ 
	gives
	\[
	0=x\cdot v_{\lambda}=\lambda(\varphi_{0})x_0\cdot v_{\lambda}
	\]
	Therefore, $x^0_0=0=x^-_0$. We shall prove that, for any $X_{\beta,i}\in\B$, 
	$x^0_{\beta, i}=0=x^-_{\beta,i}$. Proceeding by induction, 
	we assume that $x_{\gamma,j}=0=x^0_{\gamma, j}$ for any $X_{\gamma, j}\in\B$ 
	such that $\hgt{\gamma}<n$. Fix $\beta\in\Qp$ with $\hgt{\beta}=n$. 
	Then, for generic $\lambda\in\h^*$, we have
	$X^*_{\beta,i}v_{\lambda}\in M(\lambda)_{\beta}$ and, since $X_{\beta,j} X_{\beta,i}^* v_{\lambda}= \drc{ij}v_{\lambda}$, 
	\begin{align}
		0=x\cdot X^*_{\beta,i}v_{\lambda}=\sum_{j} x^-_{\beta, j}x^0_{\beta,j}X_{\beta,j}X^*_{\beta,i}v_{\lambda}=\lambda(x^0_{\beta,i})x^-_{\beta, i}v_{\lambda}
	\end{align} 
	Therefore, $x^0_{\beta,i}=0=x^-_{\beta, i}$.
\end{proof}


\subsection{Regularity of the matrix coefficients on $M(\lambda)$}\label{ss:quantum-regularity}
For any $\lambda\in\h^*$, let $M^*(\lambda)$ be the (restricted) dual Verma module
and $(\cdot,\cdot)_{M(\lambda)}:M(\lambda)\ten M^*(\lambda)\to\hext{\IC}$ the natural pairing. 
\begin{proposition}
	For any $\lambda\in\h^*$, $v\in M(\lambda)$, and $f\in M(\lambda)^*$,
	the matrix coefficient $(xv,f)_{M(\lambda)}$ lies in $\hext{\IC[\lambda]}$.
\end{proposition}

\begin{proof}
	Note that, for any $x^\pm\in \Uhnpm{}$, the coefficient
	$(x^-v, x^+f)\in\hext{\IC}$ is independent of $\lambda$.
	We can write $x=\sum_i x^+_i x^0_i x^-_i$, for some $x^+_i\in \Uhnp{}$, $x^0_i\in \hext{U\h}$, and $x^-_i\in (\Uhnm{})_{\beta_i}$, with $\beta_i\in\Qp$. 
	Then, we have
	\[
	(xv,f)_{M(\lambda)}=\sum_i(x^0_i x^-_iv,S(x_i^+)f)_{M(\lambda)}
	=\sum_i(\lambda-\beta_i)(x^0_i)(x^-_iv,S(x^+_i)f)_{M(\lambda)}\,.
	\]
	The result follows.
\end{proof}


\subsection{Proof of Theorem~\ref{ss:quantum-faithful}}\label{ss:proof-quantum-faithfulness}

Assume that $x\in\DrAh{}$ acts trivially on every category $\O$ integrable $\Uhg$--module.
We shall prove that $x$ acts trivially on any Verma module, so that
 $x=0$ by Proposition~\ref{ss:quantum-verma}.

Clearly, $x$ acts trivially on $M(\lambda)$ if and only if, for any $v\in M(\lambda)$ 
and $f\in M(\lambda)^*$, the matrix coefficient $(x v,f)_{M(\lambda)}$ vanishes.
By Proposition~\ref{ss:quantum-regularity}, it is enough to check that this holds for
$\lambda$ in a Zariski open subset of $\h^*$. To this end, note that, if $v\in M(\lambda)_{\beta}$, then $xv=x(\beta)v$, where $x(\beta)\in U\g$ is the truncation
of $x$ at $\beta$. Therefore, it is possible to choose $\lambda\in\Pp$ large enough
such that
\[
(xv,f)_{M(\lambda)}=(xv,f)_{L(\lambda)}=0
\] 
\ie $(xv,f)_{M(\lambda)}$ is equal to the matrix coefficient of $x$ on the unique irreducible quotient $L(\lambda)$ of $M(\lambda)$. By assumption on $x$, 
the latter is zero, since $L(\lambda)$ is integrable for $\lambda\in\Pp$.
The result follows.



\section{Quantum Weyl group actions of pure braid groups}\label{s:pure-qWg}

\subsection{Completions}\label{ss:completions}

Let $A$ be an algebra, $\C\subset\Rep(A)$ a full subcategory, and $\End(\sff_\C)$ the
algebra of endomorphisms of the forgetful functor $\sff_\C:\C\to\vect$. By definition, an
element of $\End(\sff_\C)$ is a collection
\[\varphi=\{\varphi_V\}_{V\in\C}\in\prod_{V\in\C}\End(V)\]
which is natural, \ie such that $f\circ\varphi_V=\varphi_W\circ f$ for any $f:V\to W$ in
$\C$. The action of $A$ on any $V\in\C$ yields a morphism of algebras $A\to\End(\sff
_\C)$, and factors through the action of $\End(\sff_\C)$ on $V$. 
We shall refer to $\End(\sff_\C)$ as the completion of $A$ with respect to 
the category $\C$. 

\subsection{Braid groups and completions}\label{ss:braid-completions}
The braid group actions considered in Section~\ref{s:km} can be 
concisely described in terms of completions.
For instance, let $\CQUOinfint$ be 
the algebra of endomorphisms 
of the forgetful functor $\FF{\hbar}\int:\OinfintUhg\to\Vect_{\hbar}$.
The quantum Weyl group
operators $\qWSk{i}$ defined by \eqref{eq:qwg-op} are elements
of $\CQUOinfintx$, and yield a group homomorphism $\QBw:\Bw
\to\CQUOinfintx$.

\subsection{Sign character of the pure braid group}\label{ss:sign}

Let $Z$ be the free abelian group with a generator $p_\alpha$ for each positive
real root $\alpha$, endowed with the $W$--action given by $w\,p_\alpha=p_{|w
\alpha|}$, where $|w\alpha|=\pm w\alpha$ according to whether $w\alpha$ is
positive or negative.

Let $\Pw\subset\Bw$ be the pure braid group. Its abelianisation $\Pwab=\Pw/[\Pw,\Pw]$
is acted upon by $\Bw/\Pw\simeq W$. By~\cite[Thm. 2.5]{tits-66} and \cite{digne-15}
the assignment $p_{\alpha_i}\to S_i^2$ uniquely extends 
to a $W$--equivariant isomorphism $Z\to\Pwab$. 

Define the sign character of $\Pw$ to be the morphism
\begin{equation}\label{eq:signh body}
\vepsh:\Pwab\to\CQUOinfintx
\qquad\qquad
\vepsh(p_\alpha)=\exp(\iota\pi h_\alpha)
\end{equation}
where $\exp(\iota\pi h_\alpha)$ is the operator acting on a weight
space of (integral) weight $\lambda$ as multiplication by $\exp(\iota\pi\lambda
(h_\alpha))$. 

\subsection{Canonical lift of the sign character}\label{ss:sign-lift}

Let $\FFp{\hbar}:\OinfUhg\to\Vect_{\hbar}$ be the forgetful functor, and consider
the morphism $\Aut(\FFp{\hbar})\to\CQUOinfintx$ corresponding to the inclusion
$\OinfintUhg\subset \OinfUhgp$. The sign character $\vepsh$ has a canonical lift
\[\Pwab\to\Aut(\FFp{\hbar})
\qquad\qquad
p_\alpha\to\exp(\iota\pi h_\alpha)
\]
which is well--defined since for any $\V\in\OinfUhg$ and $n\geq 0$, $\V/\hbar^m\V$
is a locally finite $\h$--module. We denote this lift by the same symbol.

\subsection{Pure braid group action on category $\Oinf$}\label{ss:pure-qg}

The following is one of the main results of this paper.

\begin{theorem}\label{thm:pure-qWg-action}
Let $\QBw:\Bw\to\CQUOinfint$ be the quantum Weyl group action of the braid group $\Bw$.
Then, the following holds.
\begin{enumerate}\itemsep0.1cm
\item For any $p\in\Pw$,
\[\QBw(p)=\vepsh(p)\cdot\QPwD(p)\]
where $\vepsh(p)$ is the sign character \eqref{eq:signh body}, and $\QPwD(p)$
is a unique element of $\Uhg$ which is invertible and of weight zero.
\item The assignment $p\to \QPwD(p)$ is a homomorphism $\Pw\to\left(\Uhg\right)^\h$
which is $\Bw$--equivariant.
\item 
The quantum Weyl group action of the pure braid group $\Pw$ on integrable modules
extends to an action 
\begin{equation}\label{eq:lambda eps body}
\QPw:\Pw\to\Aut(\FFp{\hbar})\qquad\text{given by}\qquad
\QPw(p)={\veps}_\hbar(p)\cdot\QPwD(p)
\end{equation}
%
\item The map $\QPw$ intertwines the inner action of $\Pw$ on $\Uhg$ \ie for any element
$Y\in\Uhg$ and $p\in\Pw$
\[\QPw(p)Y\QPw(p)^{-1}=p(Y)\]
in $\CQUOinf$.
\end{enumerate}
\end{theorem}

\begin{proof}
(2),(3) and (4) follow from (1).

(1) It suffices to prove the existence of $\QPwD(p)$ for a set of generators of $\Pw$. The
uniqueness of $\QPwD(p)$ for any $p\in\Pw$ then follows from Theorem~\ref{thm:faithful}. 
By \cite[Cor. 6]{digne-gomi-01} (see also \cite[Prop.~2.5]{digne-15}), $\Pw$ is generated
by the elements $\topS{w}\topS{i}^2\topS{w}^{-1}$, where $i\in\bfI$, $w\in W$ is such that
$w\rt{i}>0$, and $\topS{w}\in\Bw$ is the canonical lift of $w$.

Consider first the case $w=1$. By \cite[Sec.~5.2]{lusztig-93}, the square of the operator
$\qWSk{i}$ is  related to the quantum Casimir operator of $U_{\hbar}\sl{2}^{\alpha_i}=
\langle E_i, F_i, \crt{i}\rangle\subset U_{\hbar}\g$ as follows. Let $\FF{\hbar,i}\int:\OinfintUhsli
\to\Vecth$ be the forgetful functor. An element of $\End(\FF{\hbar,i}\int)$ is determined
by its action on each of the indecomposable representations $\{\V_r^i\}_{r\geq 0}$, where $\V_r^i$ is of rank $r+1$. The
Casimir operator $\Cu{i}$ of $U_{\hbar}\sl{2}^{\alpha_i}$
	acts on $\V_r^i$ as multiplication by $\symd{i} r(r+2)/2$. Set $\ku{i}= \Cu{i}-\symd{i}\crt{i}^2/2$, so that
	$\ku{i}$ acts on the subspace of $\V_r^i$ of weight $m\alpha_i/2$ as multiplication by
	$\symd{i}( r(r+2)-m^2)/2$.
	Then,
\begin{equation}\label{eq:Si2}
\qWSk{i}^2=\exp(\iota\pi\crt{i})\cdot q^{\tm{\ku{i}}}
\end{equation}

	By \cite[Sec.~5]{drinfeld-89}, 
	\begin{equation}\label{eq:Drinfeld Cas}
	q^{\Cu{i}}=\sum_{m\geqslant0}F_i^m\phi_mE_i^m
	\end{equation}
	for some explicit $\phi_m\in\hext{U\h_i}$. It follows that $q^{\Cu{i}}$
	lies in $\Uhg$, 
	and therefore so does $q^{\ku{i}}=q^{\Cu{i}}q^{-\symd{i}\crt{i}^2/2}$.
	Thus, setting $\QPwD(\topS{i}^2)=q^{\tm{\ku{i}}}\in\DrAh$, we get
	\[
	\QBw(\topS{i}^2)=\qWSk{i}^2=\exp(\iota\pi\crt{i})\cdot q^{\tm{\ku{i}}}=\vepsh(\topS{i}^2)\cdot \QPwD(\topS{i}^2)
	\]
	
Note next that if $w\in W$ satisfies $w\rt{i}>0$, then $\qWT{w}=\Ad(\qWSk{w})$ satisfies
	$\qWT{w}(E_i)\in\Uhbp{w\rt{i}}$, and
	$\qWT{w}(F_i)\in\Uhbm{-w\rt{i}}$ \cite[Sec.~37.1]{lusztig-93}. It follows that
	$q^{\ku{w,i}}=\qWT{w}(q^{\ku{i}})$ is a weight zero element in $\DrAh$, and
	if we set $\QPwD(\topS{w}\topS{i}^2\topS{w}^{-1})= q^{\tm{\ku{w,i}}}$, then
	\begin{align}
		\QBw(\topS{w}\topS{i}^2\topS{w}^{-1})=\qWSk{w}\qWSk{i}^2\qWSk{w}^{-1}=\exp(\iota\pi\crt{w,i})\cdot q^{\tm{\ku{w,i}}}=\vepsh(p)\cdot \QPwD(p)
	\end{align}
\end{proof}

\begin{remarks}\label{rk:P Ug}
\begin{enumerate}
\item The proof of Theorem \ref{thm:pure-qWg-action} shows that the action $\lambda$ on category
$\Oinf$ modules for $\Uhg$ is explicitly given on the generators of $\Pw$ by
\[\QBw(S_w S_i^2S_w^{-1})=\exp(\iota\pi\crt{w,i})\cdot \tm{q^{\ku{w,i}}}
\]
\item Since $\QPwD$ maps to $\Uhg$, it defines a (signless quantum Weyl group) action of $\Pw$
on {\it any} $\Uhg$--module.
\end{enumerate}
\end{remarks}

\subsection{The normally ordered quantum Weyl group action}
\label{ss:norm-qW}

We shall be interested in the following modification of the action \eqref
{eq:lambda eps body}. Let
\[\bw:\Pw\to\exp(\hbar\h)\subset\Uhg 
\qquad\text{be given by}\qquad
\bw(p_\alpha)=q^{\tinv{\alpha}}=\exp(\hbar\tinv{\alpha}/2)\]
(cf.~Section~\ref{ss:no-qWg}). Define the morphism
\[{\noQPwO}:\Pw\to\Uhg
\qquad\text{by}\qquad
{\noQPwO}(p)=\QPwD(p)\cdot \bw(p)^{-1}\]
so that
$\QPw(p)={\veps}_\hbar(p)\cdot \noQPwO(p)\cdot \bw(p)$ for any $p\in\Pw$.

We refer to $\noQPwO$ as the {\it normally ordered} quantum Weyl group
action of $\Pw$. The terminology is justified by 
the fact that, for any $i\in\bfI$, $\noQPwO(S_i^2)$ acts 
as the normally ordered quantum Casimir operator, 
in contrast with \eqref{eq:Si2}. Namely, one has
\[
\noQPwO(S_i^2)=\QPwD(S_i^2)\cdot\bw(p_{\alpha_i})^{-1}=q^{2\Ku{i}{+}}
\]
where $\Ku{i}{+}=(\ku{i}-\tinv{\alpha_i})/2$.
This modified action will be relevant in Theorem~\ref{thm:main}.
Note also that for any 
element $Y\in\Uhg$ of weight $\gamma\in\sfQ$ and $p\in\Pw$, one has
\[\Ad(\noQPwO(p))(Y)=p(Y)\cdot({\vepsh(p)},\gamma)^{-1}\cdot(\bw(p),\gamma)^{-1}\]
in $\CQUOinf$.

\subsection{Pure braid group actions for $\Uqg$.}

Let $\bbK$ be a field of characteristic zero, $q\in\bbK^\times$ an element of infinite order, \eg
$q\in\bbC^\times$ not a root of unity or $q\in\bbQ(q)$, and $\Uqg$ the corresponding quantum
group over $\bbK$.

The definition of (integrable) category $\Oinf$ $\Uqg$--modules is similar to the formal case (see
\eg~\cite[Ch.~3]{lusztig-93}). The analogues of Theorem~\ref{thm:pure-qWg-action} and Section
\ref{ss:norm-qW} hold for $\Uqg$ and defines actions of $\Pw$ on category $\Oinf$ modules.

In this case, the quantum Casimirs $q^{\ku{i}}$ do not lie in $\Uqg$, but in the Drinfeld algebra
$\DrAq$ of $\Uqg$, and the morphism $\QPwD$ takes values in $\DrAq$. Note that the latter acts on
any category $\Oinf$ module $\V$ since, for any $v\in\V$, $(U_q\n^+)_\beta v=0$ for all but finitely
many $\beta\in\sfQ_+$.


\section{The Casimir connection}\label{s:Casimir}

\subsection{Fundamental group of root system arrangements}\label{ss:fund-group}

Let $\GCM{A}$ be a symmetrisable generalised Cartan matrix, $(\h_{\IR},\Pi,\Pi^\vee)$
a realisation of $\GCM{A}$ over $\IR$, and $(\h=\IC\ten_{\IR}\h_{\IR},\Pi,\Pi^\vee)$ its complexification.
Let $\Pi^\perp\subset\h$ be the annihilator of $\Pi$, set $\h\ess=\h/\Pi^\perp$, and note
that $\h\ess$ is independent of the realisation of $\GCM{A}$. Let
\begin{equation}
	\C=\{h\in\h\ess_{\IR}\;|\;\forall i\in\bfI,\,\rt{i}(h)>0\}
\end{equation}
be the fundamental Weyl chamber in $\h\ess_{\IR}$, and ${\sf Y}_{\IR}=\bigcup_{w\in W}w(\ol{\C})$
the 
Tits cone. ${\sf Y}_{\IR}$ is a convex cone, and the Weyl group
$W$ acts properly discontinuously on its interior $\mathring{{\sf Y}}_{\IR}$ and complexification
${\sf Y}=\mathring{{\sf Y}}_{\IR}+\iota\h\ess_{\IR}\subseteq\h\ess$ 
\cite{looijenga-80,vinberg-71}. The regular points of this action are given by
\[
\sfX={\sf Y}\setminus\bigcup_{\alpha\in\Rs{+}}\Ker(\alpha)
\]
The action of $W$ on $\sfX$ is proper and free, and the space $\sfX/W$ is a complex
manifold.
The following result is due to van der Lek \cite
{van-der-lek-83}, and generalises Brieskorn's Theorem \cite{brieskorn-71} to the case of an 
arbitrary Weyl group.

\begin{theorem}
The fundamental groups of $\sfX/W$ and $\sfX$ are isomorphic to $\Bw$ and $\Pw$ respectively.
\end{theorem}

The generators $\{S_i\}_{i\in\bfI}$ of $\Bw$ may be described as follows. Let $p:\sfX\to\sfX/W$
be the canonical projection, fix a point $x_0\in\C$ and use $[x_0]=p(x_0)$ as a base point in
$\sfX/W$. For any $i\in\bfI$, choose an open disk $D_i$ in 
$x_0+\IC\crt{i}$, centered in $x_0-\frac{\rt{i}(x_0)}{2}\crt{i}$, and 
such that $\ol{D}_i$
does not intersect any root hyperplane other than $\Ker(\alpha_i)$.
Let $\gamma_i:[0,1]\to x_0+\IC\crt{i}$ be the path from $x_0$ to $s_i(x_0)$ in $\sfX$ determined by $\gamma_i\vert_{[0,1/3]\cup[2/3,1]}$ is affine and lies in 
$x_0+\IR\crt{i}\setminus D_i$, the points $\gamma_i(1/3),\gamma_i(2/3)$ are in $\partial\ol{D}_i$, and
$\gamma_i|_{[1/3,2/3]}$ is a semicircular arc in $\partial\ol{D}_i$,
positively oriented with respect to the natural orientation of $x_0+\IC\crt{i}$. Then, $\topS{i}= p\circ\gamma_i$.

\subsection{The Casimir connection}\label{ss:casimir-conn}

For any positive root $\alpha\in\Delta_+$, let $\{e_{\pm\alpha}^{(i)}\}_{i=1}^{\rsm
{\alpha}}$ be bases of $\g_{\pm\alpha}$ which are dual with respect to $\iip{\cdot}
{\cdot}$, and
\begin{equation}
	\Ku{\alpha}{+}=\sum_{i=1}^{\rsm{\alpha}}e_{-\alpha}^{(i)}e_\alpha^{(i)}
\end{equation}
the corresponding truncated and normally ordered Casimir operator. 
Let $\V$ be a $\g$--module in category $\Ohinfg$ and
$\IV= \sfX\times\V$ the holomorphically trivial vector bundle over 
$\sfX$ with fibre $V$. Finally, set $\nablah=\frac{\hbar}{2\pi\iota}$.

\begin{definition}
	The Casimir connection of $\g$ is the connection on $\IV$ given by
	\begin{equation}\label{eq:Casimir}
		\nablak=
		d-\nablah\sum_{\alpha\in\Rs{+}}\frac{d\alpha}{\alpha}\cdot\Ku{\alpha}{+}
	\end{equation}
\end{definition}
\noindent
Note that the sum converges in the $\hbar$--adic topology since, for any $v\in\V$ and
$n\geq 0$, $\Ku{\alpha}{+} v\in\hbar^n\V$ for all but finitely many $\alpha\in\Rs{+}$.
\Omit{\footnote
{Our main interest is when $\V$ is defined over $\IC$, \ie is isomorphic as $\g$--module
to $V\fml$, where $V\in\Oinfg$. In that case, the sum over $\alpha$ is locally finite on $V$.
The constructions in this section carry over {\it verbatim} to the case of a general $\V\in
\Ohinfg$, however.}}

The Casimir connection for a semisimple Lie algebra was discovered by De
Concini around '95 (unpublished, though the connection is referenced in \cite
{procesi-96}) and, independently, Millson--Toledano Laredo \cite{toledano-02,
millson-toledano-05} and Felder--Markov--Tarasov--Varchenko \cite{felder-markov-tarasov-varchenko}.
In \cite{felder-markov-tarasov-varchenko}, the case of an arbitrary symmetrisable
Kac--Moody algebra is considered.

The connection $\nablak$ is flat (see \cite{felder-markov-tarasov-varchenko}
and \cite[Thm.~3.4]{appel-toledano-15}) and therefore yields a monodromy
representation 
\begin{align}
\PT:\Pw=\Poid(\sfX;x_0)\to\GLoid(\V)
\end{align}
Moreover, since the coefficients of $\nablak$ have weight zero, the
action of $\Pw$ preserves the generalised weight spaces of $\V$.

This is more conveniently expressed in terms of completions. Let $\sff:
\Ohinfg\to\Vect_\hbar$ be the forgetful functor. Then, the monodromy of $
\nablak$ yields an action
\begin{align}\label{eq:pure-mono-action}
	\PT:\Pw=\Poid(\sfX;x_0)\to\Aut(\sff)
\end{align}

\subsection{The orbifold fundamental groupoid of $\sfX$}
\label{ss:rep-groupoid}

Let $\Poid(\sfX;Wx_0)$ be the fundamental groupoid of $\sfX$ based at the $W
$--orbit of $x_0$. Then, $\Poid{(\sfX/W;[x_0])}$ is equivalent to the orbifold
fundamental groupoid $\WPoid{(\sfX; Wx_0)}$, which is defined as follows.
\vskip .2cm

\begin{itemize}\itemsep0.1cm
	\item Its set of objects is $Wx_0$.
	\item A morphism between $x,y\in Wx_0$ is a pair $(w,\gamma)$,
	where $w\in W$ and $\gamma$ is a path in $\sfX$ from $x$ to $w^{-1}y$.
	\item The composition of $(w,\gamma):x\to y$ and $(w',\gamma')
	: y\to z$ is given by
	\[(w',\gamma')\circ (w,\gamma)=
	(w'w, w^{-1}(\gamma')\circ\gamma):x\to z\]
\end{itemize}

The projection functor
\begin{equation}\label{eq:projection functor}
	P:\WPoid{(\sfX; Wx_0)}\longrightarrow\Poid{(\sfX/W; [x_0])}
\end{equation}
given by $P(wx_0)=[x_0]$ and $P(w,\gamma)=[\gamma]$ is fully faithful since, for any
given $x,y\in Wx_0$, a loop $
[\gamma]\in \Poid{(\sfX/W; [x_0])}$ lifts uniquely to a path $\gamma: x\to w^{-1}y$, for
a unique $w\in W$. Any $x\in Wx_0$ therefore
determines a right inverse $\cE_{x}$ of $P$ given by $\cE_{x}([x_0])=x$ and $\cE_{x}([\gamma])
=(w,\gamma)$, where $\gamma$ is the lift of $[\gamma]$ through $x$, and $w$
is such that $\gamma(1)=w^{-1}x$.

\subsection{Obstruction to $W$--equivariance \cite[Sec.~4]{appel-toledano-15}}\label{ss:equiv-obstr}

Extend the monodromy of $\nablak$ to $\Poid(\sfX;Wx_0)$, and lift it to a map
$\PT:\Poid(\sfX;Wx_0)\to\T_\g$, where $\T_\g$ is the holonomy algebra of the
root arrangement of $\g$. The lack of $W$--equivariance of $\nablak$ can then
be described by the 1--cocycle
\[\aw:W\to\Hom_{\grpd}(\Poid{(\sfX;Wx_0)},\ttg)\]
defined by
$\aw_w(\gamma)=\PT(\gamma)^{-1}\cdot w^{-1}\PT(w\gamma)$.

The following summarises the main properties of $\aw$. 
\begin{theorem}\label{th:aw}\hfill
\begin{enumerate}
\item $\aw$ is abelian, that is takes values in $\sfM=\Hom_{\grpd}(\Poid{(\sfX;Wx_0)},\exp(\hbar\h))$.
\item $\aw$ is a coboundary, that is $\aw_w=d\bw_w=\bw\cdot(w^{-1}\bw)^{-1}$ for some $\bw\in\sfM$, and any $w\in W$.
\item The cochain $\bw$ can be normalised so that $\bw(\gamma_i)=\exp(\hbar a_i\tinv{\alpha_i})$
for any given choice of $\{a_i\}_{i\in\bfI}\subset\IC$, and is then unique.
\end{enumerate}
\end{theorem}

\begin{remark}\label{rks:equi}
(1) follows from the fact that $w^{-1}\PT(w\gamma)$ is the parallel transport of
\[w^*\nablak=\nablak-\sfh a_w
\qquad\text{where}\qquad
a_w = 
\sum_{\substack{\alpha\in\Rs{+}:\\w\alpha\in\Rs{-}}}
\frac{d\alpha}{\alpha}\cdot\tinv{\alpha}\]
Since $\nablak$ and the $\h$--valued 1--form $a_w$ commute, $\aw_w$ is
the parallel transport of $d-\sfh a_w$, and in particular takes values in $\sfM$.
\end{remark}

\subsection{Equivariant monodromy \cite[Sec.~4]{appel-toledano-15}}\label{ss:equiv-conn}

For any $b\in\Bw$, let $\texpb{b}\in\CUOhinfintx$ be its action by the triple
exponentials \eqref{eq:triple}, and $\wt{b}\in\Poid(\sfX;Wx_0)$ the unique lift of
$b$ through $x_0$. The following is a direct consequence of Theorem \ref{th:aw}

\begin{theorem}\label{thm:atl15-equivariant}
There is a unique morphism $\bw:\Poid(\sfX;Wx_0)\longrightarrow \exp(\hbar\h)$
such that
\begin{enumerate}\itemsep0.25cm
\item The assignment
\[\MBwtb:\Bw\to\CUOhinfintx\qquad\qquad 
\MBwtb(b)=\texpb{b}\cdot\noMPwO(\wt b)\cdot\bw(\wt b)\]
is a group homomorphism.
\item For any $i\in\bfI$, $\bw(\gamma_i)=\exp(\hbar\tinv{\alpha_i}/4)$.
\end{enumerate}
\end{theorem}

\begin{remarks}
\begin{itemize}[leftmargin=2em]\itemsep0.25cm
\item The normalisation of $\bw(\gamma_i)$ is chosen so that, if $\g=\sl{2}$ with simple
root $\alpha_i$,
\begin{equation}\label{eq:normalisation}
\MBwtb(S_i)=
\wt{s}_i\cdot\exp(\hbar\Ku{\alpha_i}{+}/2)\cdot\exp(\hbar\tinv{\alpha_i}/4)=
\wt{s}_i\cdot\exp(\hbar\Ku{\alpha_i}{}/4)
\end{equation}
where $\Ku{\alpha_i}{}=e_if_i+f_ie_i$ is the truncated Casimir of $\sl{2}$.
\item
We shall refer to $\MBwtb$ as the \emph{monodromy action} of $\Bw$. This is justified by the
fact that, when $\g$ is of finite or affine type, $\bw$ is the monodromy of the connection $d-\sfh
A$, where $A$ is a resummation of the formal abelian 1--form
		\[\wh{A}=\frac{1}{2}\sum_{\alpha\in\Rs{+}}
		\frac{d\alpha}{\alpha}\cdot m_\alpha\tinv{\alpha}\]
		(cf.~\cite[Prop.~4.9 and Appendix A]{appel-toledano-15}). Thus, in these cases,
		$\MBwtb$ is the monodromy of the pushdown of the connection $\nablak-\sfh A$ 
		to the quotient $\sfX/W$.
\end{itemize}
\end{remarks}

\subsection{Monodromy action of the pure braid group on category $\Oinf$}\label{ss:pure-monodromy}

Let 
\begin{equation}\label{eq:sign-character}
\veps:\Pwab\to\CUOhinfintx
\qquad\qquad
\veps(p_\alpha)=\exp(\iota\pi h_{\alpha})
\end{equation}
be the sign character (cf. \ref{ss:sign}), $\FFp{}:\Ohinfg\to\Vect_{\hbar}$ the forgetful functor,
and lift $\veps$ to a morphism $\Pwab\to\Aut(\FFp{})$ as in \ref{ss:sign-lift}.

\begin{proposition}\label{prop:pure-monodromy} The following holds. 
	\begin{enumerate}\itemsep0.25cm
		\item For any $\alpha\in\Rs{+}\re$, $\texpb{p_{\alpha}}=\veps(p_{\alpha})$ and $\bw(p_\alpha)=\exp(\hbar\tinv{\alpha}/2)$.
		\item
		The restriction of $\MBwtb$ to $\Pw$ lifts to an action
\begin{equation}
\MPw:\Pw\to\Aut(\FFp{})
\qquad\text{given by}\qquad
\MPw(p)={\veps}(p)\cdot\noMPwO(p)\cdot \bw(p) 
\end{equation}
	\end{enumerate}	
\end{proposition}

\begin{proof}
	(1) For any $i\in\bfI$, 
	$\texpb{S_i^2}=\texp{i}^2=\exp(\iota\pi h_i)$
	so that, for any $w\in W$ such that $w\alpha_i>0$, 
	$\texpb{S_wS_i^2S_w^{-1}}=\exp(\iota\pi h_{w\alpha_i})$.
	Thus, $\texpb{p}=\veps(p)$ for any $p\in\Pw$.

	For the second identity, it is enough to verify the relation on 
	the loops $p_{w\alpha_i}= w(p_{\alpha_i})\in\Poid(\sfX; wx_0)$, where 
	$p_{\alpha_i}=s_i(\gamma_i)\circ\gamma_i$, for $i\in\bfI$, 
	and $w\in W$ is such that $w\alpha_i>0$ (cf.~Section~\ref{ss:fund-group}).
	For $w=\id$, one has
	\begin{align}
		\bw(p_{\alpha_i})=\bw(s_i(\gamma_i))\bw(\gamma_i)
		=s_i(\aw_{s_i}(\gamma_i)^{-1}\bw(\gamma_i))\bw(\gamma_i)
		=s_i(\aw_{s_i}(\gamma_i))^{-1}
	\end{align}
	where the second equality follows from $\aw=d\bw$, and the third one from $\bw(\gamma_i)\in\exp(\IC\hbar\tinv{\alpha_i})$.
	By Remark~\ref{ss:equiv-obstr}, $\aw_v$ is the parallel transport of the abelian connection
	\begin{equation}\label{eq:Av}
	d-\nablah
	\sum_{\substack{\alpha\in\Rs{+}:\\v\alpha\in\Rs{-}}}
	\frac{d\alpha}{\alpha}\cdot\tinv{\alpha}
	\end{equation}
	For $v=s_i$, this is $d-\sfh d\log\alpha_i\cdot\tinv{\alpha_i}$, so that $\aw_{s_i}(\gamma_i)
	=\exp(\hbar\tinv{\alpha_i}/2)$.
	
	For $w\neq\id$, one has 
	\[\bw(w(p_{\alpha_i}))=w(\aw_w(p_{\alpha_i})^{-1}\bw(p_{\alpha_i}))=w(\aw_w(p_{\alpha_i}))^{-1}\exp(\hbar\tinv{w\alpha_i}/2)\]
	Note that $d\alpha/\alpha$ has a non--zero residue on the hyperplane $\alpha_i=0$ only if $\alpha=\pm\alpha_i$. It follows from \eqref{eq:Av} for $v=w$, and 
	$w\alpha_i\in \Rs{+}$ that $\aw_w(p_{\alpha_i})=1$, whence the result.
	
	(2) follows from (1) and Theorem~\ref{ss:equiv-conn}.
\end{proof}



\section{Braided Coxeter categories}\label{s:mono-recap-1}

We review below the notion of {\it braided Coxeter category} introduced in \cite
{appel-toledano-19b}.
Informally speaking, such an object is a collection of braided monoidal categories
labelled by the subdiagrams of a given diagram $\dgr$ -- in the relevant examples
the Coxeter graph of $\g$. These are equipped with relative fiber functors corresponding
to the inclusions of subdiagrams and an additional combinatorial datum -- a {\it maximal
nested set} -- which labels points at infinity in the \DCP model of the Cartan subalgebra
of $\g$ \cite{deconcini-procesi-95}. The functors corresponding to the inclusion $\emptyset\subset\dgr$ additionally carry
distinguished automorphisms -- the {\it local monodromies} -- which give rise to an action
of the generalised braid group $\Bw$.

For $\Uhg$, such a structure arises on $\OinfintUhg$ from the
$R$--matrix and quantum Weyl group operators.
For the category $\Ohinfintg$, it arises from
the dynamical coupling of the KZ and Casimir connections of 
$\g$ \cite{toledano-16}. 
This is analogous to the fact that the monodromy of the KZ equations 
gives rise to a braided tensor category structure on $\Ohinfg$ \cite{drinfeld-89},
and the fact that the canonical fundamental solutions of the Casimir 
equations constructed by Cherednik and \DCP \cite{cherednik-89,
deconcini-procesi-95} give rise to a Coxeter structure on $\Ohinfintg$
\cite{toledano-08}.

\subsection{Nested sets \cite[Sec.~5]{appel-toledano-15}}  

A {\it diagram} is an undirected graph $\dgr$ with no multiple edges or loops. A {\it subdiagram}
$B\subseteq\dgr$ is a full subgraph that is, a graph consisting of a (possibly empty) subset of
vertices of $\dgr$, together with all edges of $\dgr$ joining any two elements of it.

Two subdiagrams $B_1,B_2\subseteq\dgr$ are {\it orthogonal} if they have no vertices in common,
and no two vertices $i_1\in B_1$, $i_2\in B_2$ are joined by an edge in $\dgr$. Two subdiagrams
$B_1,B_2\subseteq\dgr$ are {\it compatible} if either one contains the other or they are orthogonal.

A {\it nested set } on $\dgr$ is a collection $H$ of pairwise compatible, connected subdiagrams
of $\dgr$ which contains the empty subdiagram and the 
connected components
of $\dgr$. We denote by $\Mns{\dgr}$ the collections of maximal nested sets on $\dgr$.

More generally, if $B'\subseteq B\subseteq \dgr$ are two subdiagrams, a nested set on $B$
{\it relative to} $B'$ is a collection of pairwise compatible subdiagrams of $B$ which contains
the connected components of $B$ and $B'$, and in which every element is compatible with,
but not properly contained in any of the connected components of $B'$.  We denote by $\Mns
{B,B'}$ the collections of maximal nested sets on $B$ relative to $B'$.

\begin{remark}
	It is well--known that when $\dgr$ is a diagram of type $\sfA_{n-1}$
	\vspace{0.25cm}
	\[
	\begin{tikzpicture}
	 \node (V1) at (0,0) {$\bullet$};
	 \node (V2) at (1,0) {$\bullet$};
	 \node (VC) at (2,0) {$\cdots$};
	 \node (VN1) at (3,0) {$\bullet$};
	 \node (VN) at (4,0) {$\bullet$};
	 \draw (V1) -- (V2) -- (VC) -- (VN1) -- (VN);
	 \node at (0,0.25) {$\scriptstyle 1$};
	 \node at (1,0.25) {$\scriptstyle 2$};
	 \node at (3,0.25) {$\scriptstyle n-2$};
	 \node at (4,0.25) {$\scriptstyle n-1$};
	\end{tikzpicture}
	\]
maximal nested sets on $\dgr$ are in bijection with complete bracketings on the non--associative
monomial $x_1x_2\cdots x_n$. Specifically, for any $1\leqslant i\leqslant j \leqslant n$, the
connected subdiagram $[i,j]\subseteq\dgr$ corresponds to the brackets $x_1\cdots (x_i\cdots
x_{j+1})\cdots x_n$, and two subdiagrams $B_1,B_2\subseteq\dgr$ are compatible if and only
if the corresponding brackets are consistent. Similarly, maximal nested sets on $\dgr$
relative to a subdiagram $B\subset\dgr$ are in bijection with {\em partially complete} bracketings,
\ie complete except for the monomials $(x_i\cdots x_{j+1})$, where $[i,j]$ is a connected
component of $B$.
\end{remark}

\subsection{Braided Coxeter categories \cite[Sec.~9]{appel-toledano-15}}\label{ss:cox-cat}\label{ss:braided-cox-category}

A {\em labelling} $\ulm$ of a diagram $\dgr$ is the assignment of an element $m_{ij}\in\{2,3,\dots,
\infty\}$ to any pair $i, j$ of distinct vertices of $\dgr$ such that $m_{ij}=m_{ji}$ and $m_{ij}=2$ if
$i$ and $j$ are orthogonal.

Let $(\dgr,\ulm)$ be a labelled diagram. A braided Coxeter category $\cCox{}$ of type $(\dgr,\ulm)$
consists of the following data
\vspace{0.25cm}
\begin{itemize}[leftmargin=2em]\itemsep0.25cm
	\item {\bf Diagrammatic categories.} For any subdiagram $B\subseteq
	\dgr$, a braided monoidal category $\C_B$.
	\item {\bf Restriction functors.} For any pair of subdiagrams $B'\subseteq B$ and
	relative maximal nested set $\F\in\Mns{B,B'}$, a tensor functor $F_\F:\C_B\to\C_
	{B'}$.\footnote{note that $F_\F$ is not assumed to be braided.}
	\item {\bf Generalised associators.} For any pair of subdiagrams $B'\subseteq B$ and 
	relative maximal nested sets $\F,G\in\Mns{B,B'}$, an isomorphism of tensor functors
	$\DCPA{\G}{\F}:F_\F\Rightarrow F_\G$.
	\item {\bf Vertical joins.} For any chain of inclusions $B''\subseteq B'\subseteq B$, $\F\in\Mns{B,B'}$,
	and $\F'\in\Mns{B',B''}$, an isomorphism of tensor functors $\redasso{\F}{\F'}: F_{\F'}\circ F_{\F}
	\Rightarrow F_{\F'\cup\F}$.
	\item {\bf Local monodromies.} For any vertex $i$ of $\dgr$ with corresponding restriction functor
	$\Fi:\C_i\to\C_\emptyset$, a distinguished automorphism $\CoxS{}{}{i}\in\Aut(\Fi)$.\footnote{note
	that $\CoxS{}{}{i}$ is not assumed to be a tensor automorphism of $\Fi$.}\\
\end{itemize}

\noindent
These data are assumed to satisfy the following properties.
\vspace{0.25cm}
\begin{itemize}[leftmargin=2em]\itemsep0.25cm
	\item {\bf Normalisation.} If $\F=\{B\}$ is the unique
	element in $\Mns{B,B}$, then $F_\F=\id_{\C_B}$ with the trivial tensor structure.
	\item {\bf Transitivity.} For any $B'\subseteq B$ and $\F,\G,
	\H\in\Mns{B,B'}$, $\DCPA{\H}{\F}=\DCPA{\H}{\G}\circ\DCPA{\G}{\F}$
	as isomorphisms $F_\F\Rightarrow F_\H$. In particular, $\DCPA{\F}{\F}
	=\id_{F_\F}$ and $\Upsilon_{\G\F}=\Upsilon_{\F\G}^{-1}$.
	\item {\bf Associativity.} For any $B'''\subseteq B''\subseteq B'\subseteq B$, 
	$\F\in\Mns{B,B'}$, $\F'\in\Mns{B',B''}$, and $\F''\in\Mns{B'',B'''}$,
	\[\redasso{\F'\cup\F}{\F''}\cdot
	\redasso{\F}{\F'}=
	\redasso{\F}{\F''\cup\F'}\cdot
	\redasso{\F'}{\F''}\]
	as isomorphisms $F_{\F''}\circ F_{\F'}\circ F_{\F}\Rightarrow F_{\F''\cup\F'\cup\F}$.
	\item {\bf Vertical factorisation.} For any $B''\subseteq B'\subseteq B$,
	$\F,\G\in\Mns{B,B'}$ and $\F',\G'\in\Mns{B',B''}$, 
	\[
	\DCPA{(\G'\cup\G)}{(\F'\cup\F)}\circ\redasso{\F}{\F'}=\redasso{\G}{\G'}
	\circ
	\left(\begin{array}{l}\DCPA{\G}{\F}\\\phantom{00}\circ\\ \DCPA{\G'}{\F'}\end{array}\right)
	\]
	as isomorphisms $F_{\F'}\circ F_{\F}\Rightarrow F_{\G'}\circ F_{\G}$.
	\item {\bf Generalised braid relations.} 
	For any $B\subseteq\dgr$, $i\neq j\in B$ and maximal nested sets $\Ki,\Kj$ on $B$ such that $\{i\}\in\Ki,
	\{j\}\in\Kj$,  the following holds in $\Aut{F_{\Ki}}$
	\[
	\underbrace{\mathsf{Ad}\left(\DCPA{i}{j}\right)(\CoxS{\sfa}{}{j})\cdot \CoxS{\sfa}{}{i}
		\cdot \mathsf{Ad}\left(\DCPA{i}{j}\right)(\CoxS{\sfa}{}{j})\cdots}_{m_{ij}}=
	\underbrace{\CoxS{\sfa}{}{i}\cdot\mathsf{Ad}
		\left(\DCPA{i}{j}\right)(\CoxS{\sfa}{}{j})\cdot \CoxS{\sfa}{}{i}\cdots}_{m_{ij}}
	\]
	where $\DCPA{i}{j}=\DCPA{\Ki}{\Kj}$ and 
	$\CoxS{\sfa}{}{i}=\Ad{\redasso{\trunc{\Ki}{}{i}}{\trunc{\Ki}{i}{}}}(\CoxS{}{}{i})\in
	\Aut{F_{\Ki}}$\footnote{$\trunc{\Ki}{}{i}$ and $\trunc{\Ki}{i}{}$ denote the {truncations} of $\Ki$ at $\{i\}$.}.
	\item {\bf Coproduct identity.}
	For any $i\in D$, the following holds in 
	$\Aut\left(\Fi\ten \Fi\right)$
	\begin{equation}\label{eq:coxcoprod-cat}
		J_i^{-1}\circ 
		\Fi(c_i)\circ\Delta(\CoxS{}{}{i})\circ J_i=
		c_{\emptyset}\circ \left(\CoxS{}{}{i}\ten \CoxS{}{}{i}\right)
	\end{equation}
	where $J_i$ is the tensor structure on $\Fi$ and $c_i, c_{\emptyset}$ are
	the opposite braidings in $\C_i$ and $\C_{\emptyset}$, respectively.\footnote
	{Given a braided monoidal category with braiding $\beta$, we set $\beta^{\scs\operatorname{op}}_{X,Y}:=\beta_{Y,X}^{-1}$.
	} 
\end{itemize}

\subsection{Representations of braid groups}\label{ss:cox-rep-braid}

Let $\BDm$ be the braid group with generators $\topS{i}$, $i\in \dgr$, and relations \eqref{eq:gen-braid}
for the labelling $\ulm$. Let $\BBm\leqslant\BDm$ be the subgroup generated by $\topS{i}$ with $i\in B$.
Finally, let $\Br{n}$ be the braid group associated to the symmetric group $\SS_n$, with generators $\topT
{1}, \dots, \topT{n-1}$, and $\brac{n}$ the set of complete bracketings  on the non--commutative monomial
$x_1 x_2\cdots x_{n}$.\\

Let $\cCox{}=(\C_B, F_{\F}, \DCPA{\F}{\G}, \redasso{\F}{\F'}, \CoxS{}{i}{})$ be a braided Coxeter category. Then, there is a family of representations
\[\lambda^{\cCox{}}_{\F, b}:\BBm\times\Br{n}\to\Aut({F_{\F}^{\boxtimes n}})\]
labelled by $B\subseteq\dgr$, $\F\in\Mns{B}$, and $b\in\brac{n}$, which is uniquely determined by the conditions
\vspace{0.25cm}
\begin{itemize}\itemsep0.25cm
	\item $\lambda^{\cCox{}}_{\F, b}(\topS{i})=\Ad({\redasso{\trunc{\F}{}{i}}{\trunc{\F}{i}{}}})
	(\CoxS{}{}{i})_{1\dots n}$ if $\{i\}\in\F$
	and $\lambda^{\cCox{}}_{\G, b}=\Ad({\DCPA{\G}{\F}})_{1\dots n}\circ\lambda^{\cCox{}}_{\F, b}$.
	\item $\lambda^{\cCox{}}_{\F, b}(\topT{i})= R^\vee_{B, i,i+1}$ if $b= x_1\cdots (x_ix_{i+1})\cdots x_n$ and $\lambda^{\cCox{}}_{\F, b'}=\Ad({\Phi_{B, b'b}})\circ\lambda^{\cCox{}}_{\F, b}$, where $\Phi_B$ and $R^\vee_{B}$ are the associativity and commutativity constraints of $\C_B$.
\end{itemize}

\subsection{Equivalence of braided Coxeter categories}
\label{ss:precox mor}\label{ss:equiv-cox}

Let $\cCox{}$, $\cCox{}'$ be two braided Coxeter categories of type $(\dgr,\ulm)$.
An equivalence $\mCox{}:\cCox{}\to\cCox{}'$ is the data of 
\vspace{0.25cm}
\begin{itemize}[leftmargin=2em]\itemsep0.25cm
	\item For any $B\subseteq\dgr$, a braided tensor equivalence $H_B: \C_B\to\C'_B$ 
	\item For any $B'\subseteq B$ and $\F\in\Mns{B,B'}$, an isomorphism $\gamma_{\F}$ of tensor functors 
	\begin{equation}\label{eq:gamma}
		\begin{tikzcd}
			\C_B
			\arrow[r, "H_B"]
			\arrow[d,"F_{\F}"']
			&
			\C'_{B}	
			\arrow[d, "F'_{\F}"]
			\arrow[dl, Rightarrow, "\gamma_{\F}"]\\
			\C_{B'}
			\arrow[r, "H_{B'}"']
			&
			\C'_{B'}	
		\end{tikzcd}
	\end{equation}
\end{itemize}
These are required to preserve the generalised associators, vertical joins, and local monodromies.
\begin{itemize}\itemsep0.25cm
	\item For any $B'\subseteq B\subseteq\dgr$ and $\F,\G\in\Mns{B,B'}$,
	\[\DCPAC{}{\G}{\F}\circ\gamma_{\F}=\gamma_{\G}\circ\DCPAC{'}{\G}{\F}\]
	as isomorphisms $F'_\F\circ H_B\Rightarrow H_{B'}\circ F_\G$.
	\item For any $B''\subseteq B'\subseteq B\subseteq\dgr$, $\F\in\Mns{B,B'}$, and $\F'\in\Mns{B',B''}$, 
	\[	\gamma_{\F'\cup\F}\circ(\redasso{\F}{\F'})'=
	\redasso{\F}{\F'}\circ
	\left(\begin{array}{l}\gamma_{\F}\\\phantom{0}\circ\\ \gamma_{\F'}\end{array}\right)
	\]
	as isomorphisms $F'_{\F'}\circ F'_{\F}\circ H_B\Rightarrow H_{B'}\circ F_{\F\cup\F'}$.
	\item For any $i\in\dgr$, $\CoxS{}{}{i}\circ\gamma_
	{\emptyset i}=\gamma_{\emptyset i}\circ S_i'$ as isomorphisms
	$F'_i\circ H_i\Rightarrow H_\emptyset\circ F_i$.\\
\end{itemize}

Let $\mCox{}:\cCox{}\to\cCox{}'$ be an equivalence of braided Coxeter categories.
Then, the representations of the braid groups $\lambda^{\cCox{}}_{\F, b}$ and
$\lambda^{\cCox{}'}_{\F, b}$
are equivalent
through the natural isomorphism $\gamma_{\F}: F'_{\F}\circ H_{B}\Rightarrow F_{\F}$.


\subsection{The braided Coxeter category $\OCox{\Uhg,\Rmx,\qWSk{}}{\sint}$}\label{ss:Cox-qWg}

Let now $\sfA$ be a symmetrisable generalised Cartan matrix, $(\h,\Pi,\Pi^\vee)$ a realisation of $\sfA$, $\g$
the corresponding \KMA and $\dgr$ its Dynkin diagram with the 
standard labelling \eqref{eq:gen-braid}, thus $\BDm=\Bw$. To simplify the exposition, we assume that $\sfA$ is
of finite or affine type, and $\h$ is its minimal realisation.

For any proper subdiagram $B\subsetneq\dgr$, we denote by $\g_B\subsetneq\g$ the subalgebra generated
by $\{e_i,f_i,h_i\}_{i\in B}$, and set $\g_\dgr=\g$.
\footnote{
	Since $\sfA$ is assumed to be
	of finite or affine type, $\g_B=\g'_B$ is the Kac--Moody algebra corresponding to the 
	Cartan submatrix $\sfA_B$.
	For a general $\sfA$, the definition of $\g_B$
	and $\Uhg_B$ requires a realisation which is {\it diagrammatic} in the sense of \cite[Sect. 2.4]{appel-toledano-15}.} 
 Similarly, we denote by $\Uhg_B\subsetneq\Uhg$ the subalgebra
topologically generated by $\{E_i,F_i,h_i\}_{i\in B}$, and set $\Uhg_\dgr=\Uhg$.

Then, the braided Coxeter category $\OCox{\Uhg,\Rmx,\qWSk{}}{\sint}$ is given by the following data.
\begin{itemize}[leftmargin=2em]\itemsep0.25cm
	\item The diagrammatic category corresponding to $B\subseteq\dgr$ is the monoidal category $\Oint
	_{\infty, \Uhg_B}$, with braiding induced by the {universal $R$--matrix} $\Rmx_B$ of $\Uhg_B$.
	%
	\item For any $B'\subseteq B$ and $\F\in\Mns{B,B'}$, $F_\F$ is the restriction functor
	$\hRes_{B'B}:\Oint_{\infty, \Uhg_B}\to\Oint_{\infty, \Uhg_{B'}}$ with the trivial tensor structure.
	\item The generalised associators and vertical joins are trivial.
	\item The local monodromy corresponding to $i\in\dgr$ is the quantum Weyl group operator $\qWSk{i}\in\CQUOinfintxD{i}$.
\end{itemize}

\begin{remarks}\hfill
	\begin{enumerate}[leftmargin=2em]\itemsep0.25cm
		\item 
		The braided Coxeter structure on $\OCox{\Uhg,\Rmx,\qWSk{}}{\sint}$
		is particularly simple in that the restriction functors, the generalised 
		associators, and the vertical join do not depend upon the choice of a 
		maximal nested set $\F\in\Mns{B,B'}$, but only on the subdiagrams $B'\subseteq B$.
		\item 
		The category $\OCox{\Uhg,\Rmx,\qWSk{}}{\sint}$ gives rise
		to a single representation of the braid group $\Bw$ (independent of $\F$) which 
		is the quantum Weyl group 
		action $\rho:\Bw\to\CQUOinfintx$ from Section~\ref{ss:braid-completions}.
		\item
		Strictly speaking, for the coproduct identity \eqref{eq:coxcoprod-cat} to hold, 
		it is necessary to consider a Cartan correction of the quantum Weyl group 
		operator $\qWSk{i}$ (cf.~\cite[Sec.~17.3]{appel-toledano-15}). 
		For simplicity, we shall gloss over this technical detail and refer the reader
		to \cite{appel-toledano-15}.
	\end{enumerate}
\end{remarks}


\subsection{The braided Coxeter category  $\OCox{\g,\nabla}{\hbar,\sint}$}\label{ss:Cox-mono-Casimir}

In \cite[Sec.~16]{appel-toledano-15}, we defined a braided Coxeter category $\OCox{\g,\nabla}{\sint}$
which underlies the equivariant monodromy of the Casimir connection, together with that of the KZ
equations for all the subalgebras $\g_B\subseteq\g$. In outline, $\OCox{\g,\nabla}{\sint}$ is described
as follows.
\vspace{0.25cm}
\begin{itemize}[leftmargin=2em]\itemsep0.2cm
	\item The diagrammatic category corresponding to $B\subseteq\dgr$ is the braided
	monoidal category $\OhinfintgD{B}$, with associativity and commutativity constraints
	given by the KZ associator $\Phi^{\nabla}_{B}$ and $R$--matrix $R^{\nabla}_{B}=\exp(\hbar\Omega_B/2)$, where $\Omega_B\in\g_B\wh{\ten}\g_B$
	is the Casimir tensor of $\g_B$,
	cf.~\cite{drinfeld-90}.
	\item For any $B'\subseteq B$ and $\F\in\Mns{B,B'}$, $F_\F$ is the standard restriction
	functor $\sff_{B'B}:\OhinfintgD{B}\to\OhinfintgD{B'}$, with tensor structure given by the
	relative twists 
	$J^{\nabla}_{\F}$ constructed in \cite{toledano-16}, see also \cite[Sec.~13]{appel-toledano-15}.
	\item For any $B'\subseteq B$ and $\F,\G\in\Mns{B,B'}$, the natural isomorphism of tensor
	functors $F_\G\Rightarrow F_\F$ is 
	given by the De Concini--Procesi (relative) associator
	$\DCPA{\F}{\G}^{\nabla}$ constructed in 
	\cite{deconcini-procesi-95}, see also \cite[Sec.~8]{appel-toledano-15}.
	\item The vertical joins are trivial.
	\item The local monodromy corresponding to any $i\in\dgr$ is the operator (cf. \eqref{eq:normalisation})
\begin{equation}\label{eq:local nabla i}
S^{\nabla}
	_{i}=\texp{i}\cdot\exp(\hbar\Ku{\alpha_i}{}/4)
\end{equation}
\end{itemize}

\begin{remark}
Contrary to the local monodromies $S^{\nabla}_{i}$, the data $(\Phi^{\nabla}_{B}, R^{\nabla}_{B}, J^{\nabla}_{\F},
\DCPA{\F}{\G}^{\nabla})$ acts on category $\Oinf$ modules. By replacing the diagrammatic categories $\OhinfintgD
{B}$ with $\OhinfgD{B}$ and excluding the $S^{\nabla}_{i}$, one obtains a braided {\em pre}--Coxeter category
$\OCox{\g,\nabla}{\hbar}$ \cite[Sec.~15]{appel-toledano-15}.
\end{remark}

In \ref{sss:mono-data-Casimir}--\ref{sss:mono-data-KZ-Casimir}, we briefly outline the construction of the relative
De Concini--Procesi associators $\DCPA{\F}{\G}^{\nabla}$ and the relative twists $J_\F^\nabla$.

%


\subsection{Monodromy data of the Casimir connection}
\label{sss:mono-data-Casimir}

Following Cherednik \cite{cherednik-89, cherednik-91} 
and De Concini--Procesi \cite{deconcini-procesi-95} (see also
\cite[Sec.~8]{appel-toledano-15}), for any
$\F\in\Mns{\dgr}$, there is a canonical
universal solution $\DCPS{\F}$ of $\nablak$ valued in $\Aut(\sff)$. 
It is uniquely determined by its prescribed 
asymptotics on a point at infinity $\sfp_\F$ corresponding to a choice
of blow--up coordinates on $\sfX$ associated to $\F$.

For any $\F,\G\in\Mns{\dgr}$, the {\it De Concini--Procesi associator} $\DCPAC{\nabla}{\F}{\G}$ is the element of $\Aut(\sff)$ defined by 
\[\DCPS{\G}(x)=\DCPS{\F}(x)\cdot\DCPAC{\nabla}{\F}{\G}\]
where $x$ lies in the fundamental Weyl chamber. 
%
The datum of the De Concini--Procesi associators yields a combinatorial
description of the equivariant monodromy of $\nablak$ as follows (cf.~\cite
[Thm.~9.3]{appel-toledano-15}). Let $\CoxS{\nabla}{}{i}$ be given by \eqref
{eq:local nabla i}. 
Then, there is a family of representations
\begin{align}
	\MBw_{\F}:\Bw\to\Aut(\sff\int)
\end{align}
labelled by $\F\in\Mns{\dgr}$, which is uniquely determined by the conditions
\begin{itemize}\itemsep0.2cm
	\item $\MBw_{\F}(\topS{i})=\CoxS{\nabla}{}{i}$ if $\{i\}\in\F$\
	\item $\MBw_{\G}=\Ad(\DCPA{\G}{\F})\circ\MBw_{\F}$
\end{itemize}
The representation $\MBw_{\F}$ is the equivariant monodromy 
of $\nablak$ computed with respect to the fundamental solution
$\DCPS{\F}$.

\subsection{Generalised associators}
\label{sss:gen associators}

For any $B\subseteq\dgr$, one similarly obtains the associators
$\DCPAC{\nabla}{\F}{\G}\in\Aut(\sff_B)$ with $\F,\G\in\Mns{B}$ which,
together with the local monodromies $\{\CoxS{\nabla}{}{i}\}_{i\in B}$,
describe the equivariant monodromy of the Casimir connection of
$\g_B$.
These associators are related to those for $\g$ as follows.
Let $\H\in\Mns{\dgr,B}$ and $\F,\G\in\Mns{B,\emptyset}$. Then,
\cite[Thm.~3.6]{deconcini-procesi-95} implies that
\begin{equation}\label{eq:forget up}
\DCPAC{\nabla}{\H\cup\G\,}{\H\cup\F}
=
\iota_{\dgr B}(\DCPAC{\nabla}{\G}{\F})
\end{equation}
where $\iota_{\dgr B}:\End(\sff_B)\to\End(\sff_\dgr)$ is 
induced by the equality $\sff_\dgr=\sff_{\dgr B}\circ\sff_{B}$.

The relative associators corresponding to an inclusion
$B'\subseteq B$ are constructed as follows. Let $\F,\G\in\Mns
{B,B'}$, choose $\H\in\Mns{B',\emptyset}$, and set
\[\DCPA{\G}{\F}^{\nabla}=
\DCPA{\G\cup\H\,}{\F\cup\H}^{\nabla}\]
One then proves that the definition is independent of the choice
of $\H$, and that $\DCPA{\G}{\F}^{\nabla}$ centralises 
$\g_{B'}$ \cite[Thm.~3.6]{deconcini-procesi-95}, and therefore
can be thought of as an automorphism of the restriction functor
$\sff_{B'B}:\OhinfgD{B}\to\OhinfgD{B'}$.

These associators satisfy the vertical factorisation since if $B''\subseteq
B'\subseteq B$, $\F,\G\in\Mns{B,B'}$, $\F',\G'\in\Mns{B',B''}$,
\[\DCPAC{\nabla}{\G\cup\G'\,}{\F\cup\F'}=
\DCPAC{\nabla}{\G\cup\G'\,}{\G\cup\F'}\cdot \DCPAC{\nabla}{\G\cup\F'\,}{\F\cup\F'}=
\iota_{BB'}(\DCPAC{\nabla}{\G'}{\F'})\cdot \DCPAC{\nabla}{\G}{\F}
\]
where the second equality follows from \eqref{eq:forget up} and the definition
of $\DCPAC{\nabla}{\G}{\F}$.

\subsection{Monodromy data of the joint KZ-Casimir system}
\label{sss:mono-data-KZ-Casimir}

The tensor structures $\{J^{\nabla}_{\F}\}_{\F\in\Mns{\dgr}}$ on the forgetful
functor $\sff=\sff_\dgr$ are obtained from the {\it dynamical} KZ equations in
$n=2$ points
\begin{equation}\label{eq:intro DKZ}
d-\left(\sfh\frac{\Omega}{z}+\mu^{(1)}\right)dz
\end{equation}
where $z=z_1-z_2$, $\mu\in\h$ and $\mu^{(1)}=\mu\otimes 1$ as follows.

These admit a canonical solution $G_0$ which is asymptotic to $z^{\sfh\Omega}$
near $z=0$. If $\mu$ is regular and real, they also admit two canonical solutions
$G_{\pm}$ which are asymptotic to $z^{\sfh\Omega_0}\cdot\exp(z\mu^{(1)})$
as $z\to\infty$ with $\Im z\gtrless 0$, where $\Omega_0$ is the projection of
$\Omega$ onto $\h\otimes\h$ \cite[Sect. 6]{toledano-16}. Define the {\it differential twist} $J_\pm(\mu)$ by
\[ J_\pm(\mu) = G_0^{-1}(z)\cdot G_{\pm}(z) \]
where $\Im z\gtrless 0$.

Then, $J_\pm(\mu)$ kills the KZ associator for $\g$. As a function of $\mu\in\C$,
where $\C$ is the fundamental Weyl chamber, $J_\pm(\mu)$ is real analytic and
varies according to the Casimir equations \cite[Sect. 7]{toledano-16}
\[ d_\h J_\pm =
\frac{\sfh}{2} \sum_{\alpha\in\Delta_+} \frac{d\alpha}{\alpha}
\left(
\Delta(\Ku{\alpha}{+}) J_\pm 
- J_\pm \left(\Ku{\alpha}{+}\otimes 1+1\otimes \Ku{\alpha}{+} \right)
\right)
\]
It follows that, for any maximal nested set $\F\in\Mns{\dgr}$, the twist
\[J^\nabla_\F=\Delta(G_\F(\mu))^{-1}\cdot J_\pm(\mu)\cdot G_\F(\mu)^{\otimes 2}\]
where $G_\F(\mu)$ is the fundamental solution of the Casimir connection
corresponding to $\F$ (see \ref{sss:mono-data-Casimir}), is independent
of $\mu\in\C$, and a tensor structure on $\sff_\dgr$.

The relative twists $J^\nabla_{\F}$ corresponding to any $B'\subseteq B$ and $\F\in\Mns{B,B'}$
are obtained by relying on vertical factorisation as follows. Fix $H\in\Mns{B',
\emptyset}$, let $F^\nabla_{\F\cup\H}$ and $F^\nabla_{\H}$ be the tensor
structures on $\sff_B,\sff_{B'}$ corresponding  to $\F\cup\H$ and $\H$
respectively. Then, define $J^\nabla_{\F}$ by
\[ \sff_{B'}(J^\nabla_{\F}) = J^\nabla_{\F\cup\H}\cdot (J^\nabla_\H)^{-1}\]
More precisely, the \rhs is a collection of natural isomorphisms
\[\sff_{B'}\left(\sff_{B'B}(U)\otimes\sff_{B'B}(V)\right)\to
\sff_{B}(U\otimes V)=
\sff_{B'}\left(\sff_{B'B}(U\otimes V)\right)
\]
defined for any $U,V\in\OhinfgD{B}$. One can prove that it satisfies the
centraliser property, \ie commutes with the action of $\g_{B'}$ \cite[Sect. 8]
{toledano-16}. Since $\sff_{B'}$ is faithful, it follows that it is of the
form $\sff_{B'}(J^\nabla_{\F})$ for a unique $J^\nabla_{\F}$. Moreover, the
latter is independent of the choice of $\H$.

\section{The equivariant monodromy theorem}\label{s:mono-recap-2}

We review in this section the main result of \cite{appel-toledano-15}, which extends
that of \cite{toledano-08,toledano-16} to the case of an arbitrary symmetrisable
\KM algebra, and yields an equivalence of braided Coxeter categories $\OCox{\g, \nabla}
{\hbar,\sint}\to\OCox{\Uhg,\Rmx,\qWSk{}}{\sint}$.
Its proof relies on the Etingof--Kazhdan equivalence, 
which is briefly reviewed in \ref{ss:ek-fun}--\ref{ss:drinf-iso-1}.

\subsection{The Etingof--Kazhdan equivalence}\label{ss:ek-fun}

In \cite[Thm.~4.2]{etingof-kazhdan-08}, Etingof and Kazhdan construct an equivalence
of categories $\FEK{}:\Ohinfg\to\OinfUhg$, together with an isomorphism $\alpha$ of
functors
\begin{equation}\label{eq:diag-5}
	\begin{tikzcd}
		\Ohinfg \arrow[rr,"\FEK{}"]\arrow[dr,"\FF{}"',""{name=U}]
		&
		&\OinfUhg \arrow[dl,"\FF{\hbar}"]
		\arrow[Rightarrow, sloped, to=U, "\alpha"']\\
		&\Vecth&
	\end{tikzcd}
\end{equation}
where $\FF{}$ and $\FF{\hbar}$ are the forgetful functors.%
\footnote{More precisely, in \cite{etingof-kazhdan-08} \nEK construct an equivalence
$\FEK{}$ between the larger categories of \DY modules over the negative Borel subalgebra
$\bm{}$, and admissible \DY modules over $\Uhbm{}$ (see also \cite[6.13]{appel-toledano-18}).
 It easily follows that $\FEK{}$ restricts to an equivalence $\Ohinfg\to\OinfUhg$, since it is the identity on Drinfeld--Yetter $\h$--modules, see~\cite[Lemma 22.11]{appel-toledano-15}. By the same argument, it also restricts to an equivalence $\Ohg\to\OUhg$.}
The equivalence $\FEK{}$ is the identity on $\h$--modules
and preserves integrability \cite[Lemma 22.9]{appel-toledano-15}. 
It therefore gives rise to
a diagram of functors in which every face commutes
\begin{equation}\label{eq:EK-diag}
	\begin{tikzcd}
	\Ohinfintg \arrow[r,"\FEK{\sint}"]\arrow[d]
	\arrow[ddd, bend right=75,swap,"\FF{}^{\sint}"]
	\arrow[dd, bend right=45, crossing over, swap,"\FF{0}^{\sint}"]
	&\OinfintUhg \arrow[d]
	\arrow[dd, bend left=45, crossing over,"\FF{\hbar,0}^{\sint}"]
	\arrow[ddd, bend left=75,"\FF{\hbar}^{\sint}"]\\
	\Ohinfg \arrow[r,"\FEK{}"]\arrow[d, "\FF{0}"]
	\arrow[dd, bend right=45,swap, crossing over,"\FF{}^{}"]
	&\OinfUhg 
	\arrow[d, swap, "\FF{\hbar,0}"]
	\arrow[dd, bend left=45, crossing over,"\FF{\hbar}"]\\
	\Wh\arrow[r,equal]\arrow[d, "\FF{\h}"]
	&   \Wh \arrow[d, swap, "\FF{\h}"]\\
	\Vecth\arrow[r,equal]&\Vecth
\end{tikzcd}
\end{equation}
where the vertical arrows are restriction functors,
and the natural isomorphisms are either trivial or induced from $\alpha$.\footnote
{The categories $\Whp$, $\Ohinfg$ and $\OinfUhgp$ naturally fit within the diagram
\eqref{eq:EK-diag}, but are omitted for simplicity.}

\subsection{The Etingof--Kazhdan isomorphism}\label{ss:drinf-iso-1}

In terms of completions, the Etingof--Kazhdan equivalence $(\FEK{}, \alpha)$ gives rise 
to an isomorphism $\Psi:\End(\sffh)\to\End(\sff)$ via the composition
\begin{equation}\label{eq:Psi body}
\End(\sffh)\longrightarrow{} \End(\sffh\circ\FEK{}) \xrightarrow{} \End(\sff)
\end{equation}
where the first isomorphism is induced by $\FEK{}$, and the second is given by $\Ad(\alpha)$.
By \eqref{eq:EK-diag}, $\EK{}{}$ restricts to an isomorphism
$\EK{}{\sint}:\CQUOinfint\to\CUOhinfint$ such that
\begin{equation}
	\begin{tikzcd}
		\CQUOinfint\arrow[r, "\EK{}{\sint}"] & \CUOhinfint \\
		\CQUOinf\arrow[r, "\EK{}{}"']\arrow[u]
		\arrow[ur,phantom, "\circlearrowleft"] & \CUOhinf \arrow[u]\\
		\CUOhss\arrow[u]  \arrow[r, equal]
		\arrow[ur,phantom, "\circlearrowleft"]& \CUOhss\arrow[u]
	\end{tikzcd}
\end{equation}
where the vertical arrows are restriction 
to category $\Oinf$ and integrable modules.


\subsection{The classical Drinfeld algebra}\label{ss:drinf-class}

Let $\DrA$ be the analogue of the Drinfeld
algebra $\DrAh$ for $\hext{U\g}$ (cf.~Section~\ref{ss:quantum-faithful}).
Namely, for any $\beta\in\Qp$, let $\B_{\beta}=\{X_{\beta,p}\}$ be a basis of $U\np{\beta}$ 
and $\B=\bigsqcup_{\beta\in\Qp}\B_{\beta}$.
Set
\begin{align}
	\DrA_0=\left
	\{\sum_{X\in\B}c_{X}X\,\colon\, c_{X}\in U\bm{}\right\}
	={\prod_{\beta\in\Qp}}U\bm{}\ten U\np{\beta}\supset U\g
\end{align}
and $\DrA=\hext{\DrA_0}$. 
The algebra structure of $\hext{U\g}$ extends to one on $\DrA$ and yields a chain
of morphisms $\hext{U\g}\subset{\DrA}\to\CUOhinf$. Proceeding as in Section~\ref{s:faithful}
one shows that $\DrA$ embeds into $\CUOhinf$ and $\CUOhinfint$.

\subsection{The monodromy theorem}\label{ss:atl3-thm}
In \cite[Thm.~22.1]{appel-toledano-15} we prove the following.

\begin{theorem}\label{thm:atl}\hfill
	\begin{enumerate}[leftmargin=2em]\itemsep0.25cm
	\item 
	There is a  canonical equivalence of braided  pre--Coxeter categories (cf.~Remark~\ref{ss:Cox-mono-Casimir})
	\[\mathbf{H}_\g=(H_B,\gamma_\F):\OCox{\g,\nabla}{\hbar}\to\OCox{\Uhg, \Rmx,\qWSk{}}{}\]
	such that 
	\begin{itemize}\itemsep0.25cm
		\item for any $B\subseteq\dgr$, the equivalence $H_B$ is  
		the Etingof--Kazhdan functor
		\[\FEK{}_{\g_B}:\OhinfgD{B}\to\OinfUhgD{B}\]
		\item for any $B'\subseteq B$ and $\F\in\Mns{B,B'}$, 
		the natural isomorphism $\gamma_\F$ is induced by 
		the action of an invertible weight zero element 
		$\gau{\F}$ in the Drinfeld algebra of $\g_B$, \ie there is a commutative
		diagram of functors
		\begin{equation}\label{eq:gamma-g}
		\begin{tikzcd}[row sep=normal]
			&&\OinfUhgD{B} 
			\arrow[dddd, dotted, bend right=15, 
			"\FF{\hbar,B}^{}"'{name=V}]
			\arrow[ddr, "\FF{\hbar,B'B}^{}"]
			\arrow[dddl, Rightarrow, crossing over, "\gamma_{\F}"']
			&\\
			\OhinfgD{B} 
			\arrow[dddd, bend right=15, swap, "\FF{B}^{}"{name=U}]
			\arrow[ddr, "\FF{B'B}^{}"] 
			\arrow[urr, "\FEK{}_{\g_B}"]
			& & &\\
			&&&\OinfUhgD{B'} 
			\arrow[ddl, "\FF{\hbar,B'}^{}"]\\
			&\OhinfgD{B'} 
			\arrow[urr, crossing over, near end, "\FEK{}_{\g_{B'}}"]
			\arrow[ddl, "\FF{B'}"]
			\arrow[Rightarrow, to=U, "\gau{\F}"]
			& & \\
			&&\Vecth&\\
			\Vecth \arrow[urr, equal] & & &\\
		\end{tikzcd}
		\end{equation}
		where the unmarked back face is the identity and the two unmarked lateral faces are the isomorphisms $\alpha$ for $\g_B$ and $\g_{B'}$.
	\end{itemize}
	\item $\mathbf{H}_{\g}$ restricts to an equivalence of braided Coxeter categories
		\[\mathbf{H}^{\sint}_\g=(H_B^{\sint},\gamma_\F):\OCox{\g,\nabla}{\hbar,\sint}\to\OCox{\Uhg, \Rmx,\qWSk{}}{\sint}\]
		where $H_B^{\sint}=\FEK{\sint}_{\g_B}$.
	\item For any $\F\in\Mns{\dgr}$, the isomorphism
	\[\EK{\F}{\sint}=\Ad(\gau{\F})\circ\EK{}{\sint}:\CQUOinfint\to\CUOhinfint\]
	intertwines the quantum Weyl group and the monodromy actions
	of $\Bw$, \ie
	\begin{equation}\label{eq:diag-0}
		\begin{tikzcd}
			& \Bw \arrow[dl, "\QBw"'] \arrow[dr, "\MBw_{\F}"]  & \\
			\CQUOinfint \arrow[rr,"\EK{\F}{\sint}"']
			&\arrow[u,phantom,sloped, "\circlearrowleft"]
			& \CUOhinfint
		\end{tikzcd}
	\end{equation}
	where $\MBw_{\F}=\MBwtb^{\F}$ denotes the monodromy action of $\Bw$
	around the point at infinity in the De Concini--Procesi 
	compactification of $\sfX$ corresponding to $\F$.
	\end{enumerate}
\end{theorem}

Since the diagrammatic equivalences are fixed, the proof amounts
to constructing suitable isomorphisms \eqref{eq:gamma}. The construction is in two steps. First, we prove that
$\OCox{\Uhg,\Rmx,\qWSk{}}{\sint}$ is equivalent to a braided 
Coxeter category $\OCox{\g,\Rmx,\qWSk{}}{\hbar,\sint}$ with 
diagrammatic categories $\OhinfintgD{B}$, $B\subseteq\dgr$, and
standard restriction functors with non--trivial tensor structures.
The equivalence is given by the diagrammatic Etingof--Kazhdan 
functors, equipped with natural isomorphisms $\wt{\gamma}_\F$
whose construction is carried out in \cite{appel-toledano-18, appel-toledano-19b}. Then, relying on the rigidity result from 
\cite{appel-toledano-19a}, we prove that $\OCox{\g,\Rmx,\qWSk{}}{\hbar,\sint}$ is equivalent to
$\OCox{\g,\nabla}{\hbar,\sint}$ with diagrammatic equivalences given by the identity functors. 
Finally, we observe that, by \cite[Thm.~10.7]{appel-toledano-19b}, the resulting isomorphisms $\gamma_{\F}$ satisfy \eqref{eq:gamma-g} for weight zero elements $\gau{\F}$ 
in the Drinfeld algebra.




\section{The monodromy theorem in category $\Oinf$}\label{s:main}

In this section, we show the equivalence of the actions of $\Pw$
constructed in Sections~\ref{s:pure-qWg} and \ref{s:Casimir}.
The proof relies on the equivariant monodromy Theorem~\ref{thm:atl},
the explicit description of  the actions of $\Pw$ from Sections~\ref
{s:pure-qWg} and \ref{s:Casimir}, and the following auxiliary result.

\subsection{Isomorphism between Drinfeld algebras}\label{ss:drinf-iso-2}

We show below that the isomorphism $\Psi:\End(\sffh)\to\End(\sff)$ \eqref{eq:Psi body}
restricts to an isomorphism $\EK{}{\Dr}:\DrAh\to\D$. Our proof closely follows Etingof and
Kazhdan's argument \cite[Rem. p. 535]{etingof-kazhdan-08} for the analogous algebra
$\DrQh=\lim_\beta \Uhg/I_\beta$ (cf. Remark \ref{rk:Qh}), and completes their affirmative
answer to a question raised by Drinfeld \cite[Question 8.2]{drinfeld-problem-92}.\footnote
{The argument in \cite{etingof-kazhdan-08} is not complete since the modules $U\g/
I_\beta$ are not equicontinuous for an arbitrary \KMA $\g$, so that the \nEK equivalence
$\FEK{}$ cannot be applied to them. In particular, the existence of an isomorphism
between $\DrQh$ and its classical counterpart raised in \cite[Question 8.2]{drinfeld-problem-92}
is not settled by \cite{etingof-kazhdan-08}. Theorem \ref{prop:drin} yields such an isomorphism
for the algebra $\DrAh$.}

For any $\beta=\sum_ik_i\rt{i}\in\Qp$, define the height of $\beta$ by $\hgt{\beta}
=\sum_i k_i$. For any $n\geqslant 0$, let $J_n\subseteq \Uhg$ be the left ideal
generated by $(\Uhnp{})_{\beta}$ with $\hgt(\beta)>n$. Set $U_{\hbar}^{(n)}=
\Uhg/J_n$, and denote by $\iota_{mn}^\hbar:U_{\hbar}^{(n)}\to U_{\hbar}^{(m)}$
($m\leqslant n$) the natural morphisms. Their classical analogues $U^{(n)}$ and
$\iota_{mn}:U_{}^{(n)}\to U_{}^{(m)}$ ($m\leqslant n$) are defined similarly for
$\hext{\Ug}$.
\begin{theorem}\label{prop:drin}
	\hfill
	\begin{enumerate}
	\item
	There is a canonical isomorphism of $\Uhg$-modules
	$\DrAh{}\simeq\lim_nU_{\hbar}^{(n)}$.
	\item
	There is a canonical isomorphism of $\hext{U\g}$-modules
	$\DrA{}\simeq\lim_nU^{(n)}$.
	\item
	$\Psi$ restricts to an
	isomorphism of algebras $\Psi^\D:\DrAh{}\to\DrA{}$.
	\end{enumerate}
\end{theorem}

\begin{proof}
(1)--(2) The action of $\DrAh{}$ on the cyclic vector yields surjective morphisms $\phi_n:\DrAh{}
\to U_{\hbar}^{(n)}$ of $\Uhg$--modules such that $\iota_{mn}^\hbar\circ\phi_n=\phi_m$. The
corresponding morphism $\phi:\DrAh{}\to\lim_nU_{\hbar}^{(n)}$ is easily seen to be an isomorphism.

(3) The algebra structure of $\DrAh$ is encoded by the morphisms between the modules $U_
{\hbar}^{(n)}$. Namely, we have a natural isomorphism
\begin{align}
	\DrAh^{\operatorname{op}}\simeq
	\End_{\Uhg}\left(\lim_n U_{\hbar}^{(n)}\right)\simeq
	\lim_m\operatorname{co}\negthinspace\lim_n\Hom_{\Uhg}(U_{\hbar}^{(n)},U_{\hbar}^{(m)})
\end{align}
(see also \cite[Appendix A.1]{appel-13}). A similar results holds for $\DrA{}$.

The module $U^{(n)}$ (resp. $U_{\hbar}^{(n)}$) does not lie in $\Ohinfg$ (resp. $\OinfUhg$) since
it is free over $U\h\fml$. However, the fact that $U\np{\beta}v=0$ (resp. $\Uhnp{\beta} v=0$) for all
but finitely many $\beta\in\Qp$ for any weight vector $v\in U^{(n)}$ (resp. $v\in U_{\hbar}^{(n)}$)
implies that $U^{(n)}$ is an equicontinuous $\g$--module and therefore a \DY module over $\bm{}$,
and that $U_{\hbar}^{(n)}$ is an admissible \DY module over $\Uhbm{}$. One can therefore apply
the 
equivalence $\FEK{}$ to $U^{(n)}$, and finds that $\FEK{}(U^{(n)})=U_{\hbar}^{(n)}$ and
$\FEK{}(\iota_{mn})=\iota_{mn}^\hbar$ \cite[Thms.~4.1--4.2]{etingof-kazhdan-08}. 
This yields a collection of natural isomorphisms
\[\Hom_{\hext{U\g}}(U^{(n)},U^{(m)})\simeq\Hom_{\Uhg}(U_{\hbar}^{(n)},U_{\hbar}^{(m)})\]
and the desired isomorphism $\EK{}{\Dr}:\DrAh\to\D$.
\end{proof}

\subsection{The monodromy theorem}\label{ss:proof-pure-mono}


\begin{theorem}\label{thm:main}
	The monodromy of the normally ordered Casimir connection on a $\g
	$--module $\V\in\Ohinfg$ is canonically equivalent to the normally ordered
	quantum Weyl group action of the  pure braid group $\Pw$ on the \nEK
	quantisation $\FEK{}(\V)\in\OinfUhg$.
\end{theorem}


\begin{proof}
Let $\F\in\Mns{\dgr}$. By Theorem~\ref{thm:atl} (3), there is a weight zero element
$\gau{\F}\in\DrA{}^\times\subset\CUOhinfx$ such that $\EK{\F}{\sint}=\Ad(\gau{\F})
\circ\EK{}{\sint}$ intertwines the quantum Weyl group and the monodromy actions
of $\Bw$, cf.~\eqref{eq:diag-0}.  We claim that this yields a commutative diagram

\begin{equation}\label{eq:main-diag-bis}
	\begin{tikzcd}
		&& \Pw \arrow[dll, "\noQPw"'] \arrow[drr, "\noMPwO^{\F}"] 
		\arrow[ddl, "\noQPwO"' very near end] \arrow[ddr, "\noMPwO^{\F}" very near end] && \\
		\CQUOinfint \arrow[rrrr,crossing over, "\EK{\F}{\sint}"]&&&& \CUOhinfint\\
		&\DrAh\arrow[rr, "\EK{\F}{\Dr}"'] &
		& \DrA \arrow[dr, hookrightarrow ] &\\
		\CQUOinf\arrow[uu]\arrow[ur, hookleftarrow]\arrow[rrrr, "\EK{\F}{}"']& &
		& &\CUOhinf\arrow[uu]
	\end{tikzcd}
\end{equation}
where $\EK{\F}{\Dr}=\Ad(\gau{\F})\circ\EK{}{\Dr}$, $\noMPwO_{\F}$
denotes the normally ordered monodromy action of $\Pw$ around the point at infinity corresponding to $\F$, and every face is commutative. 
Then, the result follows from the commutativity of the back face.

We first prove the commutativity of the top face. Since $\gau{\F}\in\DrA$ is weight zero and $\FEK{\sint}:\Ohinfintg\to\OinfintUhg$ is the identity at the level of 
$\h$--modules in $\Vecth$, 
$\EK{\F}{\sint}=\Ad(\gau{\F})\circ\EK{}{\sint}$ intertwines the characters of $\Pw$ given by $\veps(p_{\alpha})=\exp(\iota\pi h_{\alpha})$,
and  $\bw(p_{\alpha})=\exp(\hbar\tinv{\alpha}/2)$. Therefore,
by Theorem~\ref{thm:pure-qWg-action} (1) and Proposition~\ref{prop:pure-monodromy} (3), we can remove $\veps$ and $\bw$, and
obtain the result.

The commutativity of the lateral faces follows from Sections~\ref{s:pure-qWg}
and \ref{s:Casimir}.
Namely, by Theorem~\ref{thm:pure-qWg-action} (2)
and Section~\ref{ss:norm-qW}, 
the normally ordered quantum Weyl group action of the pure braid group $\Pw\subset\Bw$ factors through the Drinfeld algebra $\DrAh\subset\CQUOinf$.
%
Moreover, by definition, $\noMPwO$ is the
normally ordered monodromy action of $\Pw$, which readily factors through the classical Drinfeld algebra $\DrA\subset\CUOhinf$.

The commutativity of the bottom and front faces follows from Section~\ref{ss:drinf-iso-2}.
Namely, by Theorem~\ref{thm:faithful} (and its analogue for $\hext{U\g}$), the
restriction to integrable category $\Oinf$ modules yields the embeddings 
$\DrAh\hookrightarrow\CQUOinfint$ and $\DrA\hookrightarrow\CUOhinfint$.
Since $\gau{\F}\in\DrA$, it follows from Theorem \ref{prop:drin} that $\EK{\F}{\sint}$ also restricts to an isomorphism $\EK{\F}{\D}=\Ad(\gau{\F})\circ\EK{}{\D}\colon\DrAh\to\DrA$.

Finally, since $\DrA$ embeds in $\CUOhinfint$, the commutativity of the top,
lateral, bottom, and front faces 
yields that of the diagram
\begin{equation}
	\begin{tikzcd}
		& \Pw \arrow[dl, "\noQPwO"'] \arrow[dr, "\noMPwO^{\F}"] & \\
		\DrAh \arrow[rr,"\EK{\F}{\D}"']
		&\arrow[u,phantom,sloped,"\circlearrowleft"]
		& \DrA
	\end{tikzcd}
\end{equation}
and the result follows. 
\end{proof}

\subsection{The equivariant monodromy theorem}
\label{ss:extension}

The following is a direct consequence of Theorem~\ref{thm:main}.

\begin{theorem}\label{co:equivariant}
Let $\V$ be a $\g$--module in category $\Oinf^\hbar$, $\FEK{}(\V)$ its \nEK quantisation,
\[\MPw:\Pw\to GL(\V)
\aand
\QPw:\Pw\to GL(\FEK{}(\V))\]
the equivariant monodromy of the Casimir connection given by Proposition \ref{prop:pure-monodromy},
and quantum Weyl group action given by Theorem \ref{thm:pure-qWg-action}.

Then, $\MPw$ and $\QPw$ are canonically equivalent. Specifically, for any $\F\in\Mns{\dgr}$
the following diagram is commutative
				\begin{equation}
			\begin{tikzcd}
				& \Pw \arrow[dl, swap, "\QPw"] \arrow[dr, "\MPw^{\F}"] & \\
				\End(\FFp{\hbar}) \arrow[rr,"\EK{\F}{}"']
				&\arrow[u,phantom,sloped,"\circlearrowleft"]
				& \End(\FFp{})
			\end{tikzcd}
		\end{equation}
\end{theorem}

\subsection{Extension to other Lie associators} 

Although Theorem \ref{ss:proof-pure-mono} and Corollary \ref{co:equivariant}
are formulated in terms of the tensor equivalence $\FEK{}:\hOinf\to\OinfUhg$
corresponding to the KZ associator, they hold true for the \nEK equivalence
corresponding to an arbitrary Lie associator $\sfPhi$. 

Indeed, by \cite{appel-toledano-18, appel-toledano-19b} the braided Coxeter category
$\OCox{\Uhg,\Rmx,\qWSk{}}{\sint}$ underlying the $R$--matrix and \qWg of $\Uhg$
(see \ref{ss:Cox-qWg}) is equivalent to a braided Coxeter category $\OCox{\g,\Rmx,
\qWSk{}}{\hbar,\sint,\sfPhi}$ with diagrammatic categories $\{\OhinfintgD{B}\}_{B\subseteq\dgr}$,
and standard restriction functors, with the corresponding horizontal equivalences
$\OhinfgD{B}\to\OinfUhgD{B}$ given by the \nEK tensor equivalence $\FEK{\sfPhi}
_{\g_B}$ corresponding to $\g_B$ and the choice of $\sfPhi$.

By the rigidity result of \cite{appel-toledano-19a}, $\OCox{\g,
\Rmx,\qWSk{}}{\hbar,\sint,\sfPhi}$ is equivalent to $\OCox{\g,\nabla}{\hbar,\sint}$,
with diagrammatic equivalences given by the identity functors endowed with a non-trivial
tensor structure. 

Composing yields an equivalence $\OCox{\g,\nabla}{\hbar,\sint}\to\OCox{\Uhg,\Rmx,
\qWSk{}}{\sint}$ whose diagrammatic equivalences are the \nEK functors corresponding
to $\sfPhi$, which then yields Theorem \ref{ss:proof-pure-mono} and Corollary \ref{co:equivariant}
for $\FEK{\sfPhi}$.



\nc{\AJ}{\aw\vert_{W_\bfJ}}
\nc{\AgJ}{\aw_\bfJ}

\nc{\BJ}{\bw_{W_\bfJ}}
\nc{\BgJ}{\bw_\bfJ}

\nc{\Jin}{_{[\bfJ]}}
\nc{\Jout}{^{[\bfJ]}}
\nc{\eJin}{\veps\Jin}
\nc{\eJout}{\veps\Jout}
\nc{\ehJout}{\veps_{\hbar}\Jout}
\nc{\AJin}{\aw\Jin}
\nc{\BJin}{\bw\Jin}
\nc{\BJout}{\bw\Jout}

\nc{\PJ}{\mathscr{P}_{\tauJin, \BJin}}
\nc{\QPBwJb}{\QPw_{\eJout,\BJout}}

\nc{\tauJin}{\tau\Jin}
\nc{\taugJ}{\tau_\bfJ}

\section{Parabolic pure braid group actions}\label{s:parabolic}

In this section, we extend the results of Sections~\ref{s:pure-qWg} and \ref{s:main}
to parabolic pure braid groups. 

\subsection{The group $\PBJ$}

For any subset $\bfJ\subseteq\bfI$, let $\PBJ\subseteq\Br{W}$ be the preimage of $\WJ
=\langle s_j\rangle_{j\in\bfJ}$ under the projection $\Br{W}\to W$. Thus, $\PB_\emptyset
=\Pw$ and $\PB_\bfI=\Bw$. The {\it parabolic pure braid group} $\PBJ$ is generated by the
braid group $\Br{W_{\bfJ}}$ and the pure braid group $\Pw$. Moreover, as an abstract
group,
\[\PBJ\simeq(\Pw\rtimes\Br{W_\bfJ})/\wt{\P}_{W_\bfJ}
\quad\text{where}\quad
\wt{\P}_{W_\bfJ}=\{(p,p^{-1})\,\vert\, p\in\P_{W_\bfJ}\}\subset\Pw\rtimes\Br{W_\bfJ}\]

\subsection{Quantum Weyl group action of $\PBJ$}\label{ss:qWg-pbj}

Let $\Uhg_{\bfJ}\subseteq\Uhg$ be the Hopf subalgebra generated by
$\{E_j,F_j,h_j\}_{j\in\bfJ}$, and $\qPOJ\subseteq\OinfUhg$ the full subcategory
of modules whose restriction to $\Uhg_{\bfJ}$ is integrable. We have the
inclusions
\[\OinfintUhg\subset\qPOJ\subset\OinfUhg\] 
together with the equalities $\qPO{\emptyset-}=\OinfUhg$ and $\qPO{\bfI-}=\OinfintUhg$.
	
Let $\FFJ{\hbar}{\bfJ--}:\qPOJ\to\Vect_{\hbar}$ be the forgetful functor.
We define below and in \ref{ss:no J} two actions
\[
\QPBwJ, \QPBwJb:\PBJ\to\Aut(\FFJ{\hbar}{\bfJ--})
\]
such that 
\begin{itemize}\itemsep0.25cm
\item for 
$\bfJ=\emptyset$, they recover the \qWg action $\lambda:\Pw\to\Aut(\FF
{\hbar}{})$ from Theorem~\ref{thm:pure-qWg-action} (3) and the normally
ordered \qWg action $\noQPwO:\Pw\to\Uhg$ from Section~\ref{ss:norm-qW},
respectively.
\item for $\bfJ=\bfI$, both give the \qWg action $\lambda:\Bw\to\Aut(\FFJ
{\hbar}{})$.
\end{itemize}

Let $\FFJ{\bfJ, \hbar}{}:\OinfintUhgD{\bfJ}\to\Vect_{\hbar}$ be the 
forgetful functor and $\lambda_{\bfJ}:\Br{W_\bfJ}\to\Aut(\FFJ{\bfJ, \hbar}{})$ 
the quantum Weyl group action of $\Br{W_\bfJ}$. Let $\QBwJ:\Br{W_\bfJ}\to\Aut(\FFJ{\hbar}{\bfJ--})$ be its lift through the restriction
functor $\qPOJ\to\OinfintUhgD{\bfJ}$.

\begin{theorem}\label{thm:parabolic-qWg-action}
	The following holds.
	\begin{enumerate}\itemsep0.1cm
		\item  The quantum Weyl group action of $\PBJ$ on integrable modules in category $\OinfUhg$
		has a unique extension to an action $\QPBwJ:\PBJ\to\Aut(\FFJ{\hbar}{\bfJ--})$
		such that $\QPBwJ|_{\Br{W_\bfJ}}=\QBwJ$ and $\QPBwJ|_{\Pw}$
		is the restriction of the action \eqref{eq:lambda eps body} to $\qPOJ\subset\OinfUhg$.
		\item The map $\QPBwJ$ intertwines the inner action of $\PBJ$ on $\Uhg$, \ie for any element $Y\in\Uhg$ and $b\in\PBJ$
		\[\QPBwJ(b)Y\QPBwJ(b)^{-1}=b(Y)\]
		in $\End(\FFJ{\hbar}{\bfJ--})$.
	\end{enumerate}
\end{theorem}

\begin{proof}
The uniqueness of $\QPBwJ$ follows from the fact that $\PBJ$ is generated by $\Pw$ and $\Br{W_\bfJ}$.
To prove the existence of $\QPBwJ$, it is enough to observe that on the one hand there is a commutative diagram
	\[
	\begin{tikzcd}
		\Br{W} \arrow[r, "\QBw"]& \Aut(\FFJ{\hbar}{})\\ 
		\Br{W_{\bfJ}} \arrow[r, "\QBwJ"'] \arrow[u, hookrightarrow]& \Aut(\FFJ{\hbar}{\bfJ--}) \arrow[u]
	\end{tikzcd}
	\]
	where 
	the right vertical arrow is induced by the inclusion $\OinfintUhg\subset\qPOJ$.
	On the other, by Theorem~\ref{thm:pure-qWg-action}, the quantum Weyl group action of $\Pw$ on integrable modules extends canonically to $\OinfUhg$ and therefore to $\qPOJ\subseteq\OinfUhg$, \ie there is a commutative diagram
	\[
	\begin{tikzcd}
		\Bw \arrow[r] \arrow[r, "\QBw"]	& \Aut(\FFJ{\hbar}{})\\ 
		\Pw \arrow[u,hookrightarrow]\arrow[r,"\QBw"']& \Aut(\FFJ{\hbar}{\bfJ--}) \arrow[u]
	\end{tikzcd}
	\]
	
	The actions of $\Br{W_\bfJ}$ and $\Pw$ on $\FFJ{\hbar}{}$ give rise to an action
	of $\Pw\rtimes\Br{W_\bfJ}$, since,
	for any $p\in\Pw$ and $b\in\Br{W_\bfJ}$, one has
	\begin{align}
		\QBwJ(b)\cdot \QPw(p)&=\QBwJ(b)\cdot{\veps}_\hbar(p)\cdot\QPwD(p)\\
		&=b({\veps}_\hbar(p))\cdot b(\QPwD(p))\cdot\QBwJ(b)\\
		&={\veps}_\hbar(bpb^{-1})\cdot\QPwD(bpb^{-1})\cdot\QBwJ(b)\\
		&=\QPw(bpb^{-1})\cdot\QBwJ(b)
	\end{align}
	where the third equality follows the $\Bw$--equivariance of $\QPwD$ (Theorem~\ref{thm:pure-qWg-action} (2)).
	Moreover, they coincide on 
	$\P_{W_\bfJ}=\Pw\cap\Br{W_\bfJ}$, and therefore give rise to the
	desired action $\QPBwJ:\PBJ\to\Aut(\FFJ{\hbar}{\bfJ--})$.
\end{proof}

\subsection{Normally ordered \qWg action of $\PBJ$.}\label{ss:no J}

Let $\Rs{\bfJ}\subseteq\Rs{}$ be the root subsystem generated by $\{\alpha_j\}
_{j\in\bfJ}$, and let
\begin{equation}\label{eq:J-sign-character}
\ehJout:\PBJ\to\Aut(\FFJ{\hbar}{\bfJ--})
\qquad\mbox{and}\qquad
\BJout:\PBJ\to\exp(\hbar\h)
\end{equation}
be the morphisms uniquely defined by the following conditions.
\begin{itemize}\itemsep0.25cm
	\item For any $b\in\Br{W_\bfJ}$, $\ehJout(b)=1=\BJout(b)$.
	\item For any $\alpha\in\Rs{\bfJ,+}\re$,
	$\ehJout(p_\alpha)=1=\BJout(p_\alpha)$.
	\item For any $\alpha\in\Rs{+}\re\setminus\Rs{\bfJ,+}\re$,
	$\ehJout(p_\alpha)=\exp(\iota\pi h_\alpha)$ and $\BJout(p_\alpha)=\exp(\hbar\tinv{\alpha}/2)$.
\end{itemize}

Note that $\ehJout$ and $\BJout$ are both $\Br{W_\bfJ}$--equivariant. They
therefore give rise to a morphism 
\[\QPBwJb:\PBJ\to\Aut(\FFJ{\hbar}{\bfJ--})\qquad\mbox{by}\qquad \QPBwJ(b)=\ehJout(b)\cdot \QPBwJb(b)\cdot \BJout(b)\] 
for any $b\in\PBJ$, which we shall refer to as the normally ordered \qWg action of $\PBJ$
on $\qPOJ$. If $\bfJ=\emptyset$, $\QPBwJb$ is the action of $\Pw$ constructed in \ref{ss:norm-qW}
while, if $\bfJ=\bfI$, $\QPBwJb$ is the \qWg action of $\Bw$ on $\OinfintUhg$.

\subsection{Tits extension and $\PBJ$}\label{ss:pbj-Tits}

Let now $\g_\bfJ\subseteq\g$ be the subalgebra generated by $\{e_j,f_j\}_{j\in\bfJ}$, 
$\POJ\subseteq\Ohinfg$ the full subcategory of modules whose restriction to $\g_\bfJ$
is integrable, and $\FFJ{}{\bfJ--}:\POJ\to\Vect_{\hbar}$ the forgetful functor.

Let $\eJout:\PBJ\to\Aut(\FFJ{}{\bfJ--})$ be the sign character defined as in  \eqref
{eq:J-sign-character}, and define $\eJin:\Pw\to\Aut(\FFJ{}{\bfJ--})$ by the relation 
\[
\veps(p)=\eJin(p)\cdot\eJout(p)
\]
for any $p\in\Pw$, where $\veps$ is given by \eqref{eq:sign-character}. Thus, 
$\eJin(p_\alpha)=\exp(\iota\pi h_\alpha)$ if $\alpha\in\Rs{\bfJ,+}\re$, and $\eJin
(p_\alpha)=1$ if $\alpha\notin\Rs{\bfJ,+}\re$.

\begin{lemma}\label{lem:pbj-triple}
Let $\V$ be a module in $\POJ$. 
Then, there is an action $\tauJin$ of $\PBJ$ on $\V$ uniquely determined
by the following conditions.
\begin{enumerate}\itemsep0.25cm
	\item The restriction of $\tauJin$ to $\Br{W_\bfJ}$ is given by
	the action $\taugJ$ of the triple exponentials \eqref{eq:triple}
	indexed by $\bfJ$.
	\item The restriction of $\tauJin$ to $\Pw$ is given by
	the sign character $\eJin$.
\end{enumerate}
\end{lemma}

\begin{proof}
	The result follows at once from Proposition~\ref{prop:pure-monodromy} (1).
\end{proof}

\remark 
Equivalently, $\tauJin$ is given by a projection of $\PBJ$ onto the 
Tits extension $\wt{W}_\bfJ$. Note also that, for $\bfJ=\emptyset$, $\tauJin$ is trivial, while, for $\bfJ=\bfI$, $\tauJin=\tau$.

\subsection{Monodromy action of $\PBJ$ on category $\O_{\infty}^{\scsop{\bfJ-int}}$}\label{ss:mono-pbj}

We construct below an action
\[\PJ:\PBJ\to\Aut(\FFJ{}{\bfJ--})\]
by making the monodromy of the Casimir connection $\nabla_{\K}$ of $\g$ equivariant, as
described in \ref{ss:intro equiv} and \ref{ss:equiv-obstr}--\ref{ss:pure-monodromy},
but only \wrt the parabolic
subgroup $W_\bfJ$. For $\bfJ=\emptyset$, $\PJ$ is the monodromy
$\PT:\Pw\to\Aut(\FF{}{})$ of $\nabla_{\K}$ (cf. Section~\ref{ss:casimir-conn}) while,
for $\bfJ=\bfI$, $\PJ$ is the equivariant monodromy action $\MBwtb:\Bw\to\Aut(\FFJ
{}{})$ of Theorem~\ref{thm:atl15-equivariant}. 

Let $\PT:\Poid{(\sfX;W x_0)}\to\T_\g$ be the monodromy of $\nabla_{\K}$, where
$\T_\g$ is the image of the holonomy algebra (cf.~\ref{ss:intro equiv}),
and consider its restriction to 
$\Poid{(\sfX;W_\bfJ x_0)}$.
The lack of equivariance of $\PT$ under $W_\bfJ$ is controlled by the 1--cocycle
\[\AJin=i_\bfJ^*\,\AJ:W_\bfJ\to\Hom_{\grpd}(\Poid{(\sfX;W_\bfJ x_0)},\exp(\hbar\h))\]
where $i_\bfJ:\Poid{(\sfX;W_\bfJ x_0)}\to\Poid{(\sfX;W x_0)}$ is the inclusion.

The obstruction $\AJin$ is related to the one for the Casimir connection 
of $\g_\bfJ$ as follows. Consider the quotient map 
\[p_\bfJ:\h\ess\to\h\ess/\bigcap_{\alpha\in\Rs{\bfJ}}\Ker(\alpha)\simeq\h_\bfJ\ess\]
$p_\bfJ$ is equivariant under $W_\bfJ$ and, by \cite[Prop. 3.12]{kac-90}, restricts
to a map $\sfX\to\sfX_\bfJ$ of Tits cones. It therefore induces a morphism of groupoids
$p_\bfJ:\Poid{(\sfX;W_\bfJ x_0)}\to\Poid{(\sfX_\bfJ;W_\bfJ [x_0]_\bfJ)}$, where
$[x_0]_\bfJ=p_\bfJ(x_0)$, which we denote by the same symbol.

\begin{lemma}
Let
\[\AgJ:W_\bfJ\to\Hom_{\grpd}(\Poid{(\sfX_\bfJ;W_\bfJ [x_0]_{\bfJ})},\exp(\hbar\h_\bfJ))\]
be the 1--cocycle measuring the lack of equivariance of the Casimir connection of
$\g_\bfJ$ \wrt $W_\bfJ$. Then, 
$\AJin=p_\bfJ^*\AgJ$.
\end{lemma}
\begin{proof}
Let $w\in W_\bfJ$. By Remark \ref{rks:equi}, $\aw_w$ is the monodromy of the 
connection $d-\sfh a_w$, where
\[\sfh a_w
 	=\nabla_{\K}-w^*\nabla_{\K}
	=\sfh\sum_{\substack{\alpha\in\Rs{+}:\\w\alpha<0}}\frac{d\alpha}{\alpha}\tinv{\alpha}
	=\sfh\sum_{\substack{\alpha\in\Rs{\bfJ, +}:\\w\alpha<0}}\frac{d\alpha}{\alpha}\tinv{\alpha}
	=p_\bfJ^*(\nabla_{\K,\bfJ}-w^*\nabla_{\K,\bfJ})\]
\end{proof}

By Theorem~\ref{thm:atl15-equivariant} for $\g_\bfJ$, $\AgJ=d\BgJ$, where $\BgJ\in
\Hom(\Poid{(\sfX_\bfJ;W_\bfJ [x_0]_\bfJ)},\exp(\hbar\h_\bfJ))$. Set $\BJin=p_\bfJ^*\BgJ$. 
Then,
\[\AJin=p_\bfJ^*\AgJ=p_\bfJ^*d\BgJ=dp_\bfJ^*\BgJ=d\BJin\]
It follows that $\BJin$ gives rise to a $W_\bfJ$--equivariant morphism 
\[\PT_{\BJin}:\Poid{(\sfX;W_\bfJ x_0)}\to\T_\g
\qquad\qquad
\PT_{\BJin}(\gamma)=\PT(\gamma)\cdot\BJin(\gamma)\]

Consider next the equivalence of groupoids
\[
P_\bfJ:
W_\bfJ\ltimes\Poid{(\sfX;W_\bfJ x_0)} 
\to
\Poid{(\sfX/W_\bfJ;[x_0])}
\cong
\PBJ
\]
generalising \eqref{eq:projection functor}. Composing with $P_\bfJ^{-1}$ 
yields a morphism 
$\PBJ\to W_\bfJ\ltimes\T_\g$ and its lift 
$\PBJ\to \PBJ\ltimes\T_\g$.
Combining this with the action $\tauJin$ of $\PBJ$ on $\FFJ{}{\bfJ--}$ defined in
Lemma~\ref{lem:pbj-triple}, yields the following generalisation of Theorem~\ref
{thm:atl15-equivariant}.

\begin{theorem}
	There is a morphism
	$\PJ:\PBJ\to\Aut(\FFJ{}{\bfJ--})$ given by
		\[\PJ(b)=\tauJin(b)\cdot\noMPwO(\wt b)\cdot\BJin(\wt b)\]
	where $\wt{b}\in\Poid{(\sfX; W_\bfJ x_0)}$ is the unique lift of $b$ through $x_0$.
\end{theorem}

	\noindent
	\remark
	Note that, for any $j\in\bfJ$, $\BJin(\gamma_j)=\exp(\hbar\tinv{\alpha_j}/4)$, since
	$p_{\bfJ}$ maps $\gamma_j$ to the corresponding generator of $\gamma_{\bfJ, i},
	\in\Pi_1(\sfX_\bfJ;W_\bfJ[x_0]_\bfJ)$ and $\BJin(\gamma_j)=\exp(\hbar\tinv{\alpha_j}/4)$
	by construction.

\subsection{The monodromy theorem for $\PBJ$}

\begin{theorem}\label{thm:parabolic}
	The $W_\bfJ$--equivariant monodromy of the Casimir connection on a $\g
	$--module $\V\in\POJ$ is canonically equivalent to the normally ordered
	quantum Weyl group action of the  parabolic braid group $\PBJ$ on the \nEK
	quantisation $\FEK{}(\V)\in\qPOJ$.
\end{theorem}

\begin{proof}
The result follows from the combination of Theorem~\ref{thm:atl} for $\Br{W_\bfJ}$
and Theorems~\ref{thm:main}--\ref{ss:extension} for $\Pw$.	
	
Specifically, let 	$B\subseteq\dgr$ be the subdiagram corresponding to $\bfJ$, $\F$
a maximal nested set containing $B\subseteq\dgr$ corresponding to $\bfJ$, and
$\F_{\bfJ}$ the induced maximal nested set on $B$. 
Let 
\[\FF{\bfJ}{}:\OhinfgD{\bfJ}\to\Vect_{\hbar}\aand
\FF{\bfJ, \hbar}{}:\OinfUhgD{\bfJ}\to\Vect_{\hbar}\]
be the forgetful functors. 
	{By Theorem~\ref{thm:atl} (1) and \eqref{eq:gamma-g}, the isomorphism $\EK{\F}{}$
	restricts to $\EK{\F_{\bfJ}}{}$, \ie there is a commutative diagram}
		\begin{equation}
		\begin{tikzcd}
			\End(\FF{\bfJ, \hbar}{}) \arrow[rr,"\EK{\F_{\bfJ}}{}"]\arrow[d]
			&
			{}
			& \End(\FFJ{\bfJ}{})\arrow[d]\\
			\End(\FF{\hbar}{}) \arrow[rr,"\EK{\F}{}"'] 
			&
			\arrow[u,phantom,sloped, "\circlearrowleft"] & \End(\FF{}{})
		\end{tikzcd}
	\end{equation}
	where the vertical arrows are induced by the restriction functors $\Ohinfg\to\OhinfgD{\bfJ}$
	and $\OinfUhg\to\OinfUhgD{\bfJ}$, respectively.
	
	Further, since the \nEK equivalence preserves the categories of $\O^{\scsop{\bfJ--int}}_\infty$ modules,
	$\EK{\F}{}$ restricts to an isomorphism $\EK{\F}{\scsop{\bfJ--int}}: \End(\FFJ{\hbar}{\bfJ--})\to\End(\FFJ{}{\bfJ--})$ such that
	\begin{enumerate}\itemsep0.25cm
		\item 
		There is a commutative diagram
		\begin{equation}
			\begin{tikzcd}
				\End(\FFJ{\bfJ, \hbar}{}) \arrow[rr,"\EK{\F_{\bfJ}}{\sint}"] \arrow[d]
				&
				{}
				& \End(\FFJ{\bfJ}{}) \arrow[d]\\
				\End(\FFJ{\hbar}{\bfJ--}) \arrow[rr,"\EK{\F}{\scsop{\bfJ--int}}"'] &
				\arrow[u,phantom,sloped, "\circlearrowleft"] & \End(\FFJ{}{\bfJ--})
			\end{tikzcd}
		\end{equation}
		where the vertical arrows are induced by the restriction functors
		$\POJ\to\OhinfintgD{\bfJ}$ and $\qPOJ\to\OinfintUhgD{\bfJ}$, respectively.
		\item There is a commutative diagram
		\begin{equation}
			\begin{tikzcd}
				\End(\FF{\hbar}) \arrow[rr,"\EK{\F}{}"] \arrow[d]
				&
				{}
				& \End(\FF{}) \arrow[d]\\
				\End(\FFJ{\hbar}{\bfJ--}) \arrow[rr,"\EK{\F}{\scsop{\bfJ--int}}"'] &
				\arrow[u,phantom,sloped, "\circlearrowleft"] & \End(\FFJ{}{\bfJ--})
			\end{tikzcd}
		\end{equation}
		where the vertical arrows are induced by the inclusions
		$\POJ\to\Ohinfg$ and $\qPOJ\to\OinfUhg$, respectively.
	\end{enumerate}

	We claim that $\EK{\F}{\scsop{\bfJ--int}}$ intertwines the actions of $\Br{W_\bfJ}$ and $\Pw$, and therefore that of $\PBJ$.
	To this end,  consider the diagram
	\begin{equation}\label{eq:diag-parabolic-0}
		\begin{tikzcd}
			&& \Br{W_\bfJ}\arrow[dd, hookrightarrow] \arrow[dll, "\lambda_{\bfJ}"'] \arrow[drr, "\noMPwO_{\tau_\bfJ, \BgJ}^{\F_{\bfJ}}"] 
			 && \\
			\End(\FFJ{\bfJ,\hbar}{}) \arrow[dd]\arrow[rrrr,crossing over, "\EK{\F_{\bfJ}}{\sint}\qquad\quad"]&&&& \End(\FFJ{\bfJ}{}) \arrow[dd]\\
			& &\PBJ \arrow[dll, "\QPBwJb"' near start] \arrow[drr, "\PJ^{\F}" near start]& &\\
			\End(\FFJ{\hbar}{\bfJ--})\arrow[rrrr, "\EK{\F}{\scsop{\bfJ--int}}"']& &
			& &\End(\FFJ{}{\bfJ--})
		\end{tikzcd}
	\end{equation}
	The front face commutes by (1); the top face by Theorem~\ref{thm:atl} (3) for $\g_\bfJ$; the left lateral face
	by Theorem~\ref{ss:qWg-pbj} (1). For the right lateral  face, recall that, for any $b\in\PBJ$, 
	\[\PJ^{\F}(b)=\tauJin(b)\cdot\noMPwO^{\F}(\wt b)\cdot\BJin(\wt b)\]
	Let $b\in\Br{W_\bfJ}$. By Lemma~\ref{lem:pbj-triple} (1), we have that $\tauJin(b)=\taugJ(b)$. 
	Then, by Remark~\ref{ss:mono-pbj}, $\BJin(\wt b)=\BgJ(\wt{b}_\bfJ)$, where $\wt{b}_\bfJ\in\Poid{(\sfX_\bfJ, W_\bfJ[x_0]_\bfJ)}$ denotes the unique
	lift of $b$ through $[x_0]_\bfJ$. Finally, $\noMPwO^{\F}(\wt b)=\noMPwO^{\F_\bfJ}(\wt{b}_\bfJ)$ since the monodromy in the 
	De Concini--Procesi compactification is recursive in nature \cite[Thm.~3.6]{deconcini-procesi-95}.
	Thus, $\EK{\F}{\scsop{\bfJ--int}}$ intertwines the actions of $\Br{W_\bfJ}$ through $\PBJ$.
	

	Similarly, consider the diagram
		\begin{equation}\label{eq:diag-parabolic-1}
		\begin{tikzcd}
			&& \Pw \arrow[dd, hookrightarrow] \arrow[dll, "\QPBwJb"'] \arrow[drr, "\noMPwO_{\eJin,\BJin}^{\F}"] 
			&& \\
			\End(\FF{\hbar}{}) \arrow[dd]\arrow[rrrr,crossing over, "\EK{\F}{}\qquad"]&&&& \End(\FF{}{}) \arrow[dd]\\
			& & \PBJ \arrow[dll, "\QPBwJb"'] \arrow[drr, "\PJ^{\F}"]  & &\\
			\End(\FFJ{\hbar}{\bfJ--})\arrow[rrrr, "\EK{\F}{\scsop{\bfJ--int}}"']& &
			& &\End(\FFJ{}{\bfJ--})
		\end{tikzcd}
	\end{equation}
	Let $p\in\Pw$ and recall the identities
		\[\veps(p)=\eJin(p)\cdot\eJout(p)\quad\mbox{and}\quad\bw(p)=\BJin(p)\cdot\BJout(p)\]
	from \ref{ss:pbj-Tits} and Remark~\ref{ss:mono-pbj}. The commutativity of the top face then follows 
	from Theorem~\ref{thm:main} by  correcting  simultaneously $\noQPw$ and $\PT^\F$ by $\eJin$ and $\BJin$.
	The left lateral face commutes by Theorem~\ref{ss:qWg-pbj} (1). The right lateral face commutes by Lemma~\ref{lem:pbj-triple} (2).
	Thus, $\EK{\F}{\scsop{\bfJ--int}}$ intertwines the actions of $\Pw$ through $\PBJ$.
\end{proof}

The (non normally ordered) \qWg action of $\PBJ$ admits a similar monodromic interpretation,
in analogy with Theorem~\ref{ss:extension}. Namely, define $\PJs:\PBJ\to\Aut(\FFJ{\hbar}{\bfJ--})$
by
\begin{equation}\label{eq:corr-par-mono}
\PJs(b)=\eJout(b)\cdot\PJ(b)\cdot\BJout(b)
\end{equation}
for any $b\in\PBJ$. Then, the following holds.

\begin{corollary}\label{co:parabolic-equivariant}
Let $\V$ be a $\g$--module in category $\Oinf^{\hbar, \scsop{\bfJ--int}}$, $\FEK{}(\V)$ its \nEK quantisation,
	\[\PJs:\PBJ\to GL(\V)
	\aand
	\QPw:\PBJ\to GL(\FEK{}(\V))\]
the corrected $W_\bfJ$--equivariant monodromy of the Casimir connection \eqref{eq:corr-par-mono},
and the quantum Weyl group action given by Theorem \ref{ss:qWg-pbj} respectively. Then, $\PJs$
and $\QPBwJ$ are canonically equivalent. 

	%
\end{corollary}


\providecommand{\bysame}{\leavevmode\hbox to3em{\hrulefill}\thinspace}
\providecommand{\MR}{\relax\ifhmode\unskip\space\fi MR }
\providecommand{\MRhref}[2]{%
	\href{http://www.ams.org/mathscinet-getitem?mr=#1}{#2}
}
\providecommand{\href}[2]{#2}

\end{document}